\documentclass[notitlepage,10pt]{article}

\usepackage{amssymb}
\usepackage{amsmath}
\usepackage{amscd}
\usepackage{latexsym}
\usepackage{fullpage}
\usepackage{eufrak}
\usepackage[all]{xy}

\begin{document}

\newtheorem{remark}{Remark}[section]
\newtheorem{remarks}[remark]{Remarks} 
\newtheorem{definition}[remark]{Definition} 
\newtheorem{lemma}[remark]{Lemma} 
\newtheorem{theorem}[remark]{THEOREM} 
\newtheorem{proposition}[remark]{PROPOSITION} 
\newtheorem{corollary}[remark]{Corollary} 
\newtheorem{recall}[remark]{Recall}
\newtheorem{notation}{Notation}[subsection]


\newcommand{\vqi}{(  V,   q, I)}
\newcommand{\vprqpripr}{(  V',   q', I')}
\newcommand{\vprqpri}{(  V',   q', I)}
\newcommand{\vqpri}{(  V,   q', I)}

\newcommand{\covqi}{  C_0(  V,  q, I)}
\newcommand{\covqpri}{  C_0(  V,  q', I)}
\newcommand{\covprqpripr}{  C_0(  V',  q', I')}

\newcommand{\simvqivprqpripr}{\hbox{\rm Sim}[\vqi,\vprqpripr]}
\newcommand{\isomvqivprqpripr}{\hbox{\rm Iso}[\vqi,\vprqpripr]}

\newcommand{\isomvqvprqpri}{\hbox{\rm Iso}[\vqi,\vprqpri]}
\newcommand{\simvqvqpri}{\hbox{\rm Sim}[\vqi,\vqpri]}
\newcommand{\isomvqvqpri}{\hbox{\rm Iso}[\vqi,\vqpri]}
\newcommand{\sisomvqvqpri}{\hbox{\rm S-Iso}[\vqi,\vqpri]}

\newcommand{\simvqvqi}{\hbox{\rm Sim}[\vqi,\vqi]}
\newcommand{\isomvqvqi}{\hbox{\rm Iso}[\vqi,\vqi]}
\newcommand{\sisomvqvqi}{\hbox{\rm S-Iso}[\vqi,\vqi]}

\newcommand{\isomcovqicovprqpripr}{\hbox{\rm Iso}[\covqi,\covprqpripr]}

\newcommand{\isomcovqcovqpri}{\hbox{\rm Iso}[\covqi,\covqpri]}
\newcommand{\isomcovqcovqprprimei}{\hbox{\rm Iso}'[\covqi,\covqpri]}
\newcommand{\sisomcovqcovqpri}{\hbox{\rm S-Iso}[\covqi,\covqpri]}

\newcommand{\sovqi}{{\hbox{\rm SO}}\vqi} 
\newcommand{\ovqi}{{\hbox{\rm O}\vqi}} 
\newcommand{\govqi}{{\hbox{\rm GO}\vqi}}

\newcommand{\autvqi}{{\hbox{\rm Aut}}\vqi} 
\newcommand{\autcovqi}{{\hbox{\rm Aut}}(\covqi)}   

\newcommand{\sautcovqi}{{\hbox{\rm S-Aut}}(\covqi)}   
\newcommand{\autcovqprimei}{{\hbox{\rm Aut}'}(\covqi)}   

\newcommand{\isolambdatwovivpripr}{{\hbox{\rm Iso}}[{\Lambda}^2(V)\otimes I^{-1},{\Lambda}^2(V')\otimes {(I')}^{-1}]}
\newcommand{\autlambdatwovi}{\hbox{\rm Aut}(\Lambda^2(V)\otimes I^{-1})}
\newcommand{\autdetlambdatwovi}{\hbox{\rm Aut}[\hbox{\rm det}(\Lambda^2(V)\otimes I^{-1})]}

\newcommand{\falgW}{\hbox{{\sf Alg}}_{{W}}}
\newcommand{\algW}{\hbox{\rm Alg}_{{W}}}
\newcommand{\faaWst}{\hbox{{\sf Id-Assoc}}_{{{W}}}}
\newcommand{\aaWst}{\hbox{{\rm Id-Assoc}}_{{{W}}}}
\newcommand{\faaWw}{{\hbox{\sf Id-}w\hbox{\sf -Assoc}_{{{W}}}}}
\newcommand{\aaWw}{{\hbox{\rm Id-}w\hbox{\rm -Assoc}_{{{W}}}}}
\newcommand{\stabw}{\hbox{\rm Stab}_w}
\newcommand{\gl}[1]{\hbox{\rm GL}_{#1}}
\newcommand{\fazuWw}{\hbox{\sf Id-}w\hbox{\sf -Azu}_{{W}}}
\newcommand{\azuWw}{\hbox{\rm Id-}w\hbox{\rm -Azu}_{{W}}}
\newcommand{\azuWiwi}{\hbox{\rm Id-}{w_i}\hbox{\rm -Azu}_{{W_i}}}
\newcommand{\azuWprwpr}{\hbox{\rm Id-}{w'}\hbox{\rm -Azu}_{{W'}}}
\newcommand{\fazuWst}{\hbox{\sf Azu}_{{W}}}
\newcommand{\azuWst}{\hbox{\rm Azu}_{{W}}}
\newcommand{\spazuWw}{\hbox{\rm Id-}w\hbox{\rm -Sp-Azu}_{{W}}}
\newcommand{\spazuWiwi}{\hbox{\rm Id-}{w_i}\hbox{\rm -Sp-Azu}_{{W_i}}}
\newcommand{\spazuWprwpr}{\hbox{\rm Id-}{w'}\hbox{\rm -Sp-Azu}_{{W'}}}
\newcommand{\spazuWst}{\hbox{\rm Sp-Azu}_{{W}}}

\newcommand{\spazuww}{\hbox{\rm Id-}w\hbox{\rm-Sp-Azu}_{W}}
\newcommand{\aaww}{{\hbox{\rm Id-}w\hbox{\rm -Assoc}_{W}}}
\newcommand{\azuww}{\hbox{\rm Id-}w\hbox{\rm -Azu}_{W}}
\newcommand{\aawst}{\hbox{\rm Id-Assoc}_{W}}
\newcommand{\azuwst}{\hbox{\rm Azu}_{W}}
\newcommand{\spazuwst}{\hbox{\rm Sp-Azu}_{W}}

\newcommand{\falgw}{\hbox{{\sf Alg}}_{W}}

\newcommand{\algw}{\hbox{Alg}_{W}}

\newcommand{\fazuww}{\hbox{{\sf Id}-w-\hbox{\sf Azu}}_{W}}

\newcommand{\faawst}{\hbox{{\sf Id-Assoc}}_{W}}

\newcommand{\fazuwst}{\hbox{{\sf Azu}}_{W}}

\newcommand{\fspazuwst}{\hbox{{\sf Sp-Azu}}_{W}}

\newcommand{\fspazuWst}{\hbox{{\sf Sp-Azu}}_{{W}}}

\newcommand{\spec}[1]{\hbox{Spec\,}\left({#1}\right)}

\newcommand{\sym}[2]{\hbox{\sf Sym}_{#1}\left[\left(
{{#2}_{#1}}^{\vee}\otimes_{#1}{{#2}_{#1}}^{\vee}\otimes_{#1}{{#2}_{#1}}\right)^{\vee}\right]}

\newcommand{\closure}[1]{\langle{#1}\rangle}

\title{\bf A Limiting Version of a Theorem in Cohomology}
\author{
\begin{tabular}{c}
 {\sc T. E. Venkata Balaji}\\
 {\em  Mathematisches Institut, Georg-August-Universit\"at}\\
 {\em Bunsenstrasse 3-5, D-37073 G\"ottingen, Germany}\\
 {\tt  tevbal@uni-math.gwdg.de}\\
 \\
 \\ 
{\bf Dedicated to Professor Martin Kneser}\\
\\
\\  
\end{tabular}
}
\date{June 27, 2004}
\maketitle
\abstract{This paper completes the study started in \cite{tevb-paper2}.
 Scheme-theoretic methods are used to classify line-bundle-valued 
 rank 3 quadratic bundles. The classification is done 
  in terms of schematic specialisations 
  of rank 4 Azumaya algebra bundles in the sense of Part A, 
 \cite{tevb-paper1}. For a quadratic form $q$ on 
 a rank 3 vector bundle $V$ with values in a line bundle $I$ over a 
 scheme $X$, the degree zero subalgebra $\covqi$ of the generalised 
 Clifford algebra $\widetilde{C}(V,q,I)$ of the triple $\vqi$  
 in the sense of Bichsel-Knus \cite{Bichsel-Knus}, is seen  
 to be such a specialised algebra by results in Part A, \cite{tevb-paper1}. 
 The Witt-invariant of  $\vqi$, which may be defined as the 
 the isomorphism class (as algebra bundle) of
 $\covqi$, is shown to determine   
 $\vqi$ upto tensoring by a twisted discriminant line bundle. Further, 
 each specialised algebra arises in this way upto isomorphism, so that the 
 association $\vqi\mapsto\covqi$ induces 
 a natural bijection 
 from the set of equivalence classes of line-bundle-valued  quadratic 
 forms on rank 3 vector bundles upto tensoring by a twisted discriminant 
 bundle and the set of isomorphism classes of schematic specialisations of 
 rank 4 Azumaya bundles over $X.$  This 
 statement may be viewed as a limiting version of the natural bijection 
 involving cohomology given by 
 ${\hbox{\rm\v H}}^1_{\hbox{\rm\tiny fppf}}(X, \hbox{\sf O}_3)/{\hbox{\rm\v H}}^1_{\hbox{\rm\tiny fppf}}(X, \mu_2)\cong
 {\hbox{\rm\v H}}^1_{\hbox{\rm\tiny \'etale}}(X, \hbox{\sf PGL}_2).$
 The special, usual and the 
 general orthogonal groups of $\vqi$ are 
 computed and canonically 
 determined in terms of $\autcovqi$, and it is shown that the 
 general orthogonal group is always a 
 semidirect product. Any element of $\autcovqi$ can be lifted to a 
 self-similarity, and in fact to an element of the orthogonal group provided
 the determinant of the automorphism is a square. The special orthogonal 
 group and the group of determinant 1 automorphisms of $\covqi$ 
 are naturally isomorphic. 
 A specialised algebra bundle $A$ 
 arises as $\covqi$ with $I\cong {\cal O}_X$ iff 
 $\hbox{\rm det}(A)\in 2.\hbox{\rm Pic}(X)$; and arises 
 with $q$ induced from a 
 global bilinear form iff the line subbundle generated by $1_A$ is 
 a direct summand of $A.$ 
 The use of the nice 
  technical notion of semiregularity introduced by
 Kneser in \cite{Kneser} allows working with 
  an arbitrary $X$, some (or even all)
 of whose points may have residue fields of
  characteristic two.\\[2mm]
 }

\tableofcontents


\section{\label{sec1}Overview of the Main Results}   

\begin{theorem}\label{bijectivity}
Over each scheme $X$ 
 the natural map $\vqi\mapsto\covqi$ that associates to a quadratic bundle 
 the degree zero subalgebra of its generalised Clifford algebra,  
 induces a bijection from the set of orbits 
\begin{displaymath}\mathbf{ 
\left\{ \begin{array}{c} \textrm{isomorphism classes }
(V,q,I), \textrm{ with}\\
 \textrm{rank}_X(V)=3, 
\textrm{rank}_X(I)=1\textrm{ and}\\
 q \textrm{ a quadratic form on the bundle $V$}\\
\textrm{with values in the line bundle 
$I$}\end{array}
\right\} \textrm{modulo the group} 
\left\{\begin{array}{c}
\textrm{ isomorphism classes } (L,h,J), \\
\textrm{ with  rank}_X(L)=\textrm{rank}_X(J)=1\\
\textrm{ and $h:L\otimes L\cong J$ a linear}\\
\textrm{ isomorphism; product induced by $\otimes$ } \end{array} \right\}}
\end{displaymath} to the set of algebra-bundle-isomorphism classes of 
associative unital algebra structures $A$, on vector bundles 
 of rank 4, that are Zariski-locally isomorphic to 
 even Clifford algebras of rank 3 quadratic 
 bundles. This bijection maps the subset 
\begin{displaymath}\mathbf{ 
\left\{ \begin{array}{c} \textrm{isomorphism classes }
(V,q,I), \textrm{ with}\\
 q \textrm{ semiregular}\end{array}
\right\} \textrm{modulo the group} 
\left\{\begin{array}{c}
\textrm{ isomorphism classes } (L,h,J), \\
\textrm{ as above} \end{array} \right\}}
\end{displaymath} surjectively onto
 the subset of isomorphism classes of 
  rank 4 Azumaya algebra bundles on $X.$
\end{theorem}

It was shown in Part A, \cite{tevb-paper1} that algebra bundles such as 
 $A$ in the statement above are precisely the scheme-theoretic specialisations 
 (or limits) of rank 4 Azumaya algebra bundles on $X.$ The main result there 
was the smoothness of the schematic closure of Azumaya algebra
 structures on a fixed 
 vector bundle of rank 4 over any scheme. 
 Part B  of \cite{tevb-paper1}  
 had applied this result 
 to obtain desingularisations (with good specialisation properties) 
 of certain moduli spaces over fairly general 
 base schemes, and applications to the study of degenerations of quadratic 
 forms on rank 3 vector bundles was initiated in \cite{tevb-paper2}. 
 Though \cite{tevb-paper2} only treats  forms with 
 values in the structure sheaf, in the following we  
  consider  
 line-bundle-valued quadratic forms. 
 Such a viewpoint is affordable thanks to the construction of
 Bichsel-Knus \cite{Bichsel-Knus} of the generalised Clifford algebra
 $\widetilde{C}(V,q,I).$  
 This is an $\mathbb{Z}$-graded algebra and we let $\covqi$ denote  
 its degree zero subalgebra. The other difference with \cite{tevb-paper2} 
 is to consider 
 $\covqi$ instead of the usual even Clifford algebra of a quadratic form with 
 values in the structure sheaf. We shall see later ((b1), Theorem \ref{structure-of-specialisation}) that it is necessary to consider line-bundle-valued 
 quadratic forms to obtain the surjectivity of Theorem \ref{bijectivity} in 
 those cases for which $\hbox{\rm det}(A)\not\in 2.\hbox{\rm Pic}(X).$ 

\begin{remark}\rm Assume for simplicity that $X$ is an affine scheme. 
 Following Sec.3, Chap.V, of the book of 
  Knus \cite{Knus}, it can be seen that the association of a semiregular 
 quadratic bundle to its even Clifford algebra induces a bijection from the 
 set of orbits 
\begin{displaymath}\mathbf{ 
\left\{ \begin{array}{c} \textrm{isomorphism classes }
(V,q,{\cal O}_X),\\
 \textrm{with rank}_X(V)=3, \textrm{ and}\\
 q \textrm{ semiregular}\end{array}
\right\} \textrm{modulo the group} 
\left\{\begin{array}{c}
\textrm{isomorphism classes } (L,h), \textrm{ with}\\
\textrm{rank}_X(L)=1\textrm{ and $h:L\otimes L\cong {\cal O}_X$}\\
\textrm{${\cal O}_X$-linear; product under $\otimes$} 
\end{array} \right\}}\end{displaymath} with
  the set of isomorphism classes of rank 4 Azumaya bundles on $X.$ 
 The group above 
 is also called the group of discriminant bundles on $X$ and is  denoted 
 by $\hbox{\rm Disc}(X).$ By Prop.3.2.2, Sec.3, Chap.III, \cite{Knus}, 
 it is naturally isomorphic to the cohomology (abelian group)  
 ${\hbox{\rm\v H}}^1_{\hbox{\rm\tiny fppf}}(X, \mu_2).$ We may thus call the 
 triples $(L,h,J)$ occurring in the statement of Theorem \ref{bijectivity} as
 twisted discriminant bundles.  
Further, 
  by Lemma 3.2.1, Sec.3, Chap.IV, \cite{Knus} the cohomology   
  ${\hbox{\rm\v H}}^1_{\hbox{\rm\tiny fppf}}(X, \hbox{\sf O}_3)$
 classifies the set of isomorphism classes of semiregular rank 3 
 quadratic bundles with values in the trivial line bundle. 
 On the other hand, 
 the set of isomorphism classes of rank 4 Azumaya algebras 
 on $X$ may be interpreted as the cohomology 
 ${\hbox{\rm\v H}}^1_{\hbox{\rm\tiny \'etale}}(X, \hbox{\sf PGL}_2)$
 (see page 145, Sec.5, Chap.III, \cite{Knus}). 
 Thus we have the bijection  
 ${\hbox{\rm\v H}}^1_{\hbox{\rm\tiny fppf}}(X, \hbox{\sf O}_3)/{\hbox{\rm\v H}}^1_{\hbox{\rm\tiny fppf}}(X, \mu_2)\cong
 {\hbox{\rm\v H}}^1_{\hbox{\rm\tiny \'etale}}(X, \hbox{\sf PGL}_2)$
and Theorem \ref{bijectivity} may be considered as the limiting version 
 of this bijection.\end{remark}

The following results lead to a proof of Theorem \ref{bijectivity}.
 Their proofs in turn follow from techniques of \cite{tevb-paper2}
 adapted to include 
 the case of line-bundle-valued forms. 
 The proof of the injectivity essentially follows from the next couple 
 of theorems, which describe the general, usual and special orthogonal 
 groups of $\vqi$ in terms of the (algebra) automorphism group 
 of $\covqi.$  

A bilinear form $b$, with values in a line bundle $I$, defined  
 on a vector bundle $V$ over the scheme $X$ 
induces an $I$-valued  quadratic form $q_b$ given on sections by 
 $x\mapsto b(x,x).$
 Further, $b$ also naturally defines an $\mathbb Z$-graded 
 linear isomorphism $\psi_b:\widetilde{C}(V,q_b, I)\cong
 \Lambda(V)\otimes L[I]$ 
  (see (2d), Theorem \ref{Bourbaki}). Here $L[I]:={\cal O}_X\oplus
\left(\bigoplus_{n> 0}(T^n(I) \oplus T^n(I^{-1}))\right)$ is the 
 Laurent-Rees algebra of  $I$, where elements of $V$ (resp. of $I$) 
 are declared to be of degree one (resp. of degree two). 
 In fact, $\widetilde{C}(V,0,I)=\Lambda(V)\otimes L[I]$. 
Since in general   
 a quadratic bundle $\vqi$ on a non-affine scheme $X$
may not be induced from a global $I$-valued bilinear 
 form, one is unable to identify the $\mathbb Z$--graded vector bundle 
 underlying its generalised Clifford algebra bundle with
 $\Lambda(V)\otimes L[I].$
 The following 
 result overcomes this problem.  

\begin{proposition}\label{transfer-to-lambda2}
Every isomorphism of algebra-bundles $\phi:\covqi\cong\covprqpripr$
is naturally associated to an 
 isomorphism of bundles $\phi_{\Lambda^2}:\Lambda^2(V)\otimes I^{-1}
\cong\Lambda^2(V')\otimes {(I')}^{-1}$ which induces a map 
$$\zeta_{\Lambda^2}:\isomcovqicovprqpripr\longrightarrow\isolambdatwovivpripr: \phi\mapsto\phi_{\Lambda^2}$$ 
where  $\isomcovqicovprqpripr$ is the set of algebra bundle isomorphisms.  
When $ V= V'$ and $I=I'$, we may thus denote  
  the subset of those $\phi$ for which
 $\hbox{\rm det}(\phi_{\Lambda^2})\in\autdetlambdatwovi\equiv
\Gamma(X, {\cal O}_X^*)$ is 
 a square by $\isomcovqcovqprprimei$ and those for which 
 $\hbox{\rm det}(\phi_{\Lambda^2})=1$ 
by the smaller subset $\sisomcovqcovqpri$. 
 Taking $ q= q'$ in these sets and replacing `` {\rm Iso}'' by ``{\rm Aut}''
 in their notations respectively defines the  
 groups $\autcovqi\supset\autcovqprimei\supset\sautcovqi.$  
\end{proposition}

 Denote by $\simvqivprqpripr$ the set of isomorphisms 
 from  $\vqi$ to $\vprqpripr.$ These are generalised similarities
  i.e., pairs  
  $(g,m)$ such that 
 $g:V\cong V'$ and $m:I\cong I'$ are  linear isomorphisms and 
 the following diagram commutes (where $q$ and $q'$ are considered 
 as morphisms of sheaves of {\em sets}): $$\begin{CD}
{\cal V} @>{g}>{\cong}> {\cal V'}\\
@V{q}VV @VV{q'}V\\
I @>{\cong}>m> I'
\end{CD}$$
 When $I=I'$, since an $m\in\hbox{\rm Aut}(I)$ may be thought of as 
 multiplication by a scalar
 $l\in\Gamma(X, {\cal O}_X^*)\equiv\hbox{\rm Aut}(I)$, we may call the 
 isomorphism $(g,m)$ as an $I$-similarity with multiplier $l.$ In such 
 a case we may as well  
 denote $(g,m)$ by the pair $(g, l)$ and we often write 
 $g:(V, q, I)\cong_{l}(V',q',I).$   Let  
 $\isomvqvprqpri$ be the subset of isometries (i.e., those pairs $(g,m)$ 
 with $m=\hbox{\rm Identity}$ or $I$-similarities with trivial 
 multipliers). 
  When $ V= V'$, the 
 subset of isometries with trivial determinant is denoted 
 $\sisomvqvqpri.$ On taking $ q= q'$
  these sets 
 naturally become subgroups of $\hbox{\rm Aut}(V)\times \Gamma(X, {\cal O}_X^*)
=\hbox{\rm GL}(V)\times \Gamma(X, {\cal O}_X^*)$
 and we define
\begin{center}\begin{tabular}{rcl}
$\simvqvqi=:\govqi$&$\supset$&$\isomvqvqi=:\ovqi$\\
 &$\supset$&$\sisomvqvqi=:\sovqi.$\end{tabular}\end{center}  
Of course, $\ovqi$ and $\sovqi$ may be thought of as subgroups of 
 $\hbox{\rm GL}(V)\equiv\hbox{\rm GL}(V)\times\{1\}.$ 

\begin{theorem}\label{lifting-of-isomorphisms}
 For $I$-valued quadratic forms $q$ and $q'$ on a rank 3 vector bundle 
  $V$ over a scheme $X$, we have the following 
 commuting diagram of natural maps of sets with the downward arrows 
 being the canonical inclusions, the horizontal arrows being surjective 
 and the top horizontal arrow being bijective:  
$$\begin{CD}
\sisomvqvqpri @>{\cong}>> \sisomcovqcovqpri \\
 @V{\hbox{\rm inj}}VV @VV{\hbox{\rm inj}}V \\
\isomvqvqpri @>{\hbox{\rm onto}}>> \isomcovqcovqprprimei \\
  @V{\hbox{\rm inj}}VV @VV{\hbox{\rm inj}}V \\
\simvqvqpri @>{\hbox{\rm onto}}>> \isomcovqcovqpri 
\end{CD}$$
With respect to the surjections of the 
 horizontal arrows in the diagram above, we further have the following (where 
 $l$  is the function that associates a similarity to its 
 multiplier, $\hbox{\rm det}(g,l):=\hbox{\rm det}(g)$ for 
 an $I$-similarity $g$ with multiplier $l$ and $\zeta_{\Lambda^2}$
 is the map of
 Prop.\ref{transfer-to-lambda2} above):
\begin{description}
\item[(a)] there is a family of sections
 $s_{2k+1}: \isomcovqcovqpri\longrightarrow\simvqvqpri$
 indexed by the integers 
 such that $l\circ s_{2k+1}
=\hbox{\rm det}^{2k+1}\circ\zeta_{\Lambda^2}$ and such 
 that $(\hbox{\rm det}^2
\circ s_{2k+1})\times(l^{-3}\circ s_{2k+1})=\hbox{\rm det}
\circ\zeta_{\Lambda^2}$;
\item[(b)] there is also a 
 section $s':\isomcovqcovqprprimei\longrightarrow \isomvqvqpri$ 
 such that $\hbox{\rm det}^2\circ s'=\hbox{\rm det}\circ\zeta_{\Lambda^2};$
\item[(c)]  there is a family of sections 
 $s_{2k+1}^{+}:\isomcovqcovqpri\longrightarrow\simvqvqpri$ indexed by 
 the integers which is 
 multiplicative when followed by the natural inclusions
 into $\hbox{\rm GL}(V)\times\Gamma(X, {\cal O}_X^*)$, i.e., 
 if $\phi_i\in\hbox{\rm Iso}[C_0(V,q_i, I), C_0(V,q_{i+1}, I)]$ 
 then $s_{2k+1}^{+}(\phi_2\circ\phi_1)=s_{2k+1}^{+}(\phi_2)\circ s_{2k+1}^{+}(\phi_1)\in\hbox{\rm GL}(V)\times\Gamma(X, {\cal O}_X^*).$
 Further, $l\circ s_{2k+1}^{+}
=\hbox{\rm det}^{2k+1}\circ\zeta_{\Lambda^2}$ and  $(\hbox{\rm det}^2
\circ s_{2k+1}^{+})\times(l^{-3}\circ s_{2k+1}^{+})=\hbox{\rm det}
\circ\zeta_{\Lambda^2}.$ 
\item[(d)] The maps $s_{2k+1}$ and 
 $s'$  above  may not be 
 multiplicative  but are
 mutliplicative upto $\mu_2(\Gamma(X, {\mathcal O}_X))$
 i.e., these maps followed  by the canonical 
 quotient map, on taking the quotient of 
 $\hbox{\rm GL}(V)\times\Gamma(X, {\cal O}_X^*)$ by
 $\mu_2(\Gamma(X, {\mathcal O}_X)).\hbox{\rm Id}_V\times\{1\}$,  
 become multiplicative. 
\end{description}
\end{theorem}

\begin{theorem}\label{lifting-of-automorphisms}
For a rank 3 quadratic bundle $\vqi$ on a scheme $X$, 
 one has the following natural 
 commutative diagram of groups with exact rows, where the downward arrows 
 are the  canonical inclusions and where 
 $l$ 
 is the function that associates to any $I$-(self)similarity its 
 multiplier:
\begin{displaymath}
\xymatrix{ 
  &  & \sovqi \ar[d]_{\hbox{\rm inj}} \ar[r]^{\cong\hspace*{1cm}}  &
 \sautcovqi \ar[d]^{\hbox{\rm inj}} & \\
 1\ar[r] & \mu_2(\Gamma(X, {\mathcal O}_X)) \ar[r]  \ar[d]_{\hbox{\rm inj}} & 
 \ovqi \ar[r] \ar[d]_{\hbox{\rm inj}} & \autcovqprimei \ar[d]^{\hbox{\rm inj}} 
 \ar[r]  & 1\\
 1\ar[r] & \Gamma(X, {\mathcal O}_X^{*}) \ar[r] & 
 \govqi \ar[r] \ar[d]_{\hbox{\rm det}^2\times l^{-3}} & \autcovqi \ar[d]^{\hbox{\rm det}} 
  \ar[r] 
  & 1\\
 & & \Gamma(X, {\mathcal O}_X^{*}) \ar@{=}[r]  & \Gamma(X, {\mathcal O}_X^{*}) & 
}
\end{displaymath}
Further, we have: 
\begin{description}
\item[(a)] There are splitting homomorphisms 
 $s_{2k+1}^{+}: \autcovqi\longrightarrow\govqi$ such that 
$l\circ s_{2k+1}^{+}=\hbox{\rm det}^{2k+1}$ and $(\hbox{\rm det}^2
\circ s_{2k+1}^{+})\times(l^{-3}\circ s_{2k+1}^{+})=\hbox{\rm det}.$ The 
 restriction of $s_{2k+1}^{+}$ to $\autcovqprimei$ does not necessarily take 
 values in $\ovqi$, but the further restriction to $\sautcovqi$ does 
 take values in $\sovqi.$ In particular, $\govqi$ is a semidirect product.
 The maps $s_{2k+1}$ and 
 $s'$ of Theorem \ref{lifting-of-isomorphisms} above (under the 
 current hypotheses) may not be 
 homomorphisms  but are homomorphisms upto $\mu_2(\Gamma(X, {\mathcal O}_X)).$ 
\item[(b)] Suppose $X$ is integral and $q\otimes\kappa(x)$ is semiregular
 at some point $x$ of 
 $X$ with residue field $\kappa(x).$
 Then any automorphism of $\covqi$ has determinant 1. Hence   
   $\autcovqi=\autcovqprimei=\sautcovqi$ and $\ovqi$ is the 
 semidirect product of $\mu_2(\Gamma(X, {\cal O}_X))$ and 
 $\sovqi.$
\end{description}
\end{theorem}
\noindent The proofs of the above results, and of the injectivity part of 
 Theorem \ref{bijectivity}, will be given in \S\ref{sec3}. 
As for the proof of the surjectivity part of  
Theorem \ref{bijectivity}, we have 

\begin{theorem}\label{structure-of-specialisation}
Let $X$ be a scheme and $ A$ a specialisation of rank 4 Azumaya algebra 
 bundles on $X.$  
 Let ${\cal O}_X. 1_A\hookrightarrow A$  be the line sub-bundle generated 
 by the  nowhere-vanishing 
 global section of $A$ corresponding to the unit
 for algebra multiplication.
\begin{description}
\item[(a)]
 There exist a rank 3 vector bundle 
 $V$ on $X$,  
 a quadratic form $q$ 
 on $V$ with values in the line bundle $I:=\hbox{\rm det}^{-1}(A)$, and 
 an isomorphism of algebra bundles 
  $A\cong \covqi.$ This gives the surjectivity in the statement of 
 Theorem \ref{bijectivity}.  Further,  the following linear 
 isomorphisms may be deduced: 
\begin{description}
\item[(1)] $\hbox{\rm det}(A)\otimes \Lambda^2({ V})\cong 
{ A}/{\cal O}_X. 1_{ A},$ from which follow:  
\item[(2)]  $\hbox{\rm det}(\Lambda^2({ V}))\cong{(\hbox{\rm det}({ A}))}^{\otimes-2};$
\item[(3)]  ${ V}\cong{({ A}/{\cal O}_X. 1_{ A})}^\vee\otimes
\hbox{\rm det}({ V})\otimes\hbox{\rm det}({ A});$
\item[(4)] $\hbox{\rm det}({ A}^\vee)\cong {(\hbox{\rm det}({ A}))}^{\otimes-3}\otimes {(\hbox{\rm det}({ V}))}^{\otimes-2}$ which 
implies that $\hbox{\rm det}({ A})\otimes\hbox{\rm det}({ A}^\vee)\in
2.\hbox{\rm Pic}(X).$ 
\end{description}
\item[(b)] There exists a quadratic bundle 
 $(V',q',I')$ such that $\covprqpripr\cong A$ 
 and with 
\begin{description}
\item[(1)] $I'={\cal O}_X$ iff $\hbox{\rm det}({ A})\in
2.\hbox{\rm Pic}(X)$;
\item[(2)] with $q'$ induced from a global $I'$-valued bilinear form 
 iff ${\cal O}_X.1_A$ is an ${\cal O}_X$-direct summand of $A$;
\item[(3)] with both $I'={\cal O}_X$ and with $q'$ induced from a 
 global bilinear form (with values in $I'$) iff 
 ${\cal O}_X.1_A$ is an ${\cal O}_X$-direct summand of $A$ and 
 $\hbox{\rm det}({ A})\in
2.\hbox{\rm Pic}(X).$
\end{description} 
\end{description}  
\end{theorem}

The next theorem is a key ingredient in the proof of part  
 (a) above.  It describes specialisations as bilinear forms under 
 certain conditions. 
 As a preparation towards its statement, we
 recall a few definitions and results  from Part A, \cite{tevb-paper1}. 
For a rank $n^2$ vector bundle $W$ 
 on a scheme $X$ and $w\in\Gamma(X,W)$ 
 a nowhere-vanishing global section, 
 recall that if  
  $\hbox{\rm Id-}w\hbox{\rm -Azu}_W$ is 
  the open $X$--subscheme of Azumaya 
 algebra structures on $ W$ with identity $ w$ then 
 its schematic image (or the scheme of specialisations or the limiting 
 scheme) in the bigger $X$--scheme 
$\hbox{\rm Id-}w\hbox{\rm -Assoc}_W$ of associative 
 $ w$-unital algebra structures on $ W$  is the $X$--scheme  
 $\hbox{\rm Id-}w\hbox{\rm -Sp-Azu}_W.$ 
By definition, the set of distinct specialised $w$-unital 
 algebra structures on 
 $W$ corresponds precisely to 
the set of global sections of this last scheme over 
 $X.$
 If 
 $\hbox{\rm Stab}_{w}\subset \hbox{\rm GL}_{W}$
 is the stabiliser subgroupscheme of $w$,  recall from Theorems 
 3.4 and 3.8, Part A, \cite{tevb-paper1}, that there exists a canonical 
 action of $\hbox{\rm Stab}_{w}$ on $\hbox{\rm Id-}w\hbox{\rm -Sp-Azu}_W$ such that the 
 natural inclusions \label{clubsuit} 
 $$(\clubsuit)\hspace*{1cm}\hbox{\rm Id-}w\hbox{\rm -Azu}_W\hookrightarrow \hbox{\rm Id-}w\hbox{\rm -Sp-Azu}_W \hookrightarrow \hbox{\rm Id-}w\hbox{\rm -Assoc-Alg}_W$$
 are all $\hbox{\rm Stab}_{w}$-equivariant.
  Now let
  $ V$ be a rank 3 vector bundle on the scheme $X$ and 
 $\hbox{\rm Bil}_{(V, I)}$
 be the associated  rank 9 vector bundle of 
 bilinear forms on $ V$ with values in the line bundle $I.$ 
 Let 
$\hbox{\rm Bil}^{sr}_{(V,I)}\hookrightarrow
 \hbox{\rm Bil}_{(V,I)}$
 correspond to the 
 open subscheme of semiregular bilinear forms. We say that a bilinear 
 form $b$ is semiregular if there is a trivialisation $\{U_i\}$ of $I$, 
 such that  
 over each open subscheme $U_i$,  the  quadratic form $q_i$ with values 
 in the trivial line bundle induced from $q_b|U_i$ is semiregular (it
  may turn out that a semiregular bilinear form may be degenerate). 
 This definition is independent of the choice of a 
 trivialisation, since $q_i$ is semiregular iff $\lambda q_i$ is semiregular 
 for every $\lambda\in\Gamma(U_i, {\cal O}_X^*).$ 
 Let $ W:=\Lambda^{even}( V, I):=\bigoplus_{n\geq 0}\Lambda^{2n}(V)\otimes I^{-\otimes n}$ and let $ w\in\Gamma(X, 
 W)$ be the nowhere-vanishing global section corresponding to the 
 unit for the natural multiplication in the twisted
  even-exterior algebra bundle $W.$
There is an obvious
  natural action of $\hbox{\rm GL}_{V}$ on $\hbox{\rm Bil}_{(V,I)}.$ 
 There is 
also a natural morphism of groupschemes  $\hbox{\rm GL}_{V}\longrightarrow
\hbox{\rm Stab}_{w}$ given on valued points
 by $g\mapsto \bigoplus_{n\geq 0}\Lambda^{2n}(g)\otimes\hbox{\rm Id}$ 
 and therefore  
 the natural inclusions marked by  $(\clubsuit)$ above are  
 $\hbox{\rm GL}_{V}$-equivariant. Finally, note that there is an obvious 
 involution $\Sigma$ on $\hbox{\rm Id-}w\hbox{\rm -Assoc-Alg}_W$
 given by $A\mapsto \hbox{\rm opposite}(A)$ which leaves the open 
 subscheme $\hbox{\rm Id-}w\hbox{\rm -Azu}_W$ invariant. 

\begin{theorem}\label{bilinear-forms-as-specialisations}\
\begin{description}
\item[(1)]
 Let $V$ be a rank 3 vector bundle on the scheme $X$,
 $ W:=\Lambda^{even}( V, I)$ and  $ w\in\Gamma(X, 
 W)$ correspond to 1 in the twisted even-exterior algebra bundle. 
 There is a natural $\hbox{\rm GL}_{V}$-equivariant 
 morphism of $X$--schemes 
$\Upsilon'=\Upsilon'_X: \hbox{\rm Bil}_{(V,I)}\longrightarrow
\hbox{\rm Id-}w\hbox{\rm -Assoc}_W$ whose 
 schematic image is precisely the scheme of specialisations 
$\hbox{\rm Id-}w\hbox{\rm -Sp-Azu}_W.$ Further 
if  $\Upsilon'$ factors canonically through 
$\Upsilon=\Upsilon_X: \hbox{\rm Bil}_{(V, I)}\longrightarrow
\hbox{\rm Id-}w\hbox{\rm -Sp-Azu}_W$, 
 then $\Upsilon$ is a $\hbox{\rm GL}_{V}$-equivariant
  isomorphism and it maps the $\hbox{\rm GL}_{V}$-stable 
 open subscheme $ \hbox{\rm Bil}^{sr}_{(V, I)}$
 isomorphically onto the $\hbox{\rm GL}_{V}$-stable open 
 subscheme $\hbox{\rm Id-}w\hbox{\rm -Azu}_W.$
\item[(2)] The involution 
 $\Sigma$ of  $\hbox{\rm Id-}w\hbox{\rm -Assoc-Alg}_W$ 
 defines a unique involution (also denoted by $\Sigma$) on 
 the scheme of specialisations $\hbox{\rm Id-}w\hbox{\rm -Sp-Azu}_W$
 leaving the open subscheme $\hbox{\rm Id-}w\hbox{\rm -Azu}_W$
 invariant, and therefore via the isomorphism $\Upsilon$, it defines 
 an involution on $\hbox{\rm Bil}_{(V,I)}.$ This involution is none other
 than the one on valued points given by $B\mapsto \hbox{\rm transpose}(-B).$  
\item[(3)]
For an $X$-scheme $T$, let $V_T$ (resp.$\thinspace W_T$, resp.$\thinspace I_T$) denote the 
 pullback of $V$ (resp.$\thinspace W$, resp.$\thinspace I$) to $T$,
  and let $w_T$ be the global section of $W_T$ induced by 
 $w.$ Then the base-changes of $\Upsilon'_X$ and $\Upsilon_X$ to $T$, 
 namely  $\Upsilon'_X\times_X T: \hbox{\rm Bil}_{(V, I)}\times_X T\longrightarrow
\hbox{\rm Id-}w\hbox{\rm -Assoc}_W\times_X T$ and 
 $\Upsilon_X\times_X T: \hbox{\rm Bil}_{(V, I)}\times_X T\cong
\hbox{\rm Id-}w\hbox{\rm -Sp-Azu}_W\times_X T$ may be  
 canonically identified with the corresponding ones over $T$ namely 
 $\Upsilon'_T: \hbox{\rm Bil}_{(V_T, I_T)}\longrightarrow
\hbox{\rm Id-}w_T\hbox{\rm -Assoc}_{W_T}$ and 
 $\Upsilon_T: \hbox{\rm Bil}_{(V_T, I_T)}\cong
\hbox{\rm Id-}w_T\hbox{\rm -Sp-Azu}_{W_T}.$
\end{description} 
\end{theorem}
\noindent The proofs of Theorems \ref{structure-of-specialisation} and 
 \ref{bilinear-forms-as-specialisations}  will be given in \S\ref{sec4}. 
Some notations and preliminaries are explained in \S\ref{sec2}.

\section{\label{sec2} Notations and Preliminaries}

This section collects together some definitions and results. It is essentially 
Section 2 of \cite{tevb-paper2} redone for forms with values in line bundles. 
We omit the proofs which follow by localisation and the corresponding 
 results of Sec.2, 
 \cite{tevb-paper2}. For the systematic treatment of the generalised 
 Clifford algebra and its properties, we refer the reader to the paper of 
Bichsel-Knus \cite{Bichsel-Knus}.

\paragraph{Quadratic and Bilinear Forms with Values in a Line Bundle.} 
 Let $V$ be a vector bundle (of constant positive rank)
 and  $I$ a line bundle on a scheme $X.$ 
 A bilinear form (resp. alternating form, resp. quadratic form) with values in
 $I$ on $V$ over an open set $U\hookrightarrow X$  is by 
 definition a section over $U$ of the vector bundle 
$\hbox{\rm Bil}_{(V, I)}$ (resp. of $\hbox{\rm Alt}^2_{(V,I)}$,
 resp. of $\hbox{\rm Quad}_{(V,I)}$), or equivalently, an 
 element of $\Gamma\left(U, \hbox{\rm Bil}_{({\cal V}, {\cal I})}
 :=(T^2_{{\cal O}_X}({\cal V}))^\vee\otimes{\cal I}\right)$ (resp. 
 of $\Gamma\left(U, \hbox{\rm Alt}^2_{({\cal V}, {\cal I})}
 :=(\Lambda^2_{{\cal O}_X}({\cal V}))^\vee\otimes{\cal I}\right)$, resp. 
 of   $\Gamma(U, \hbox{\rm Quad}_{({\cal V}, {\cal I})})$ ),
 where the sheaf 
 $\hbox{\rm Quad}_{({\cal V},{\cal I})}$---the
 (coherent locally-free) sheaf of ${\cal O}_X$-modules
 corresponding to the bundle $\hbox{\rm Quad}_{(V,I)}$ of $I$-valued 
 quadratic forms 
 on $V$---is defined by the exactness of the following sequence: 
$$0\longrightarrow \hbox{\rm Alt}^2_{({\cal V},{\cal I})} \longrightarrow \hbox{\rm Bil}_{({\cal V},{\cal I})} \longrightarrow \hbox{\rm Quad}_{({\cal V},{\cal I})}\longrightarrow 0.$$
 In terms of the corresponding (geometric) vector bundles over $X$, the above 
 translates into  the following 
 sequence of morphisms of vector bundles, with the first one a closed immersion and the second one a Zariski locally-trivial principal 
$\hbox{\rm Alt}^2_{(V, I)}$-bundle: 
$$\hbox{\rm Alt}^2_{(V,I)} \hookrightarrow \hbox{\rm Bil}_{(V,I)} \twoheadrightarrow \hbox{\rm Quad}_{(V,I)}.$$ 
Given a quadratic form $q\in\Gamma(U, \hbox{\rm Quad}_{({\cal V},{\cal I})})$,
 recall 
 that  the usual 
 `associated' symmetric 
 bilinear form $b_q\in\Gamma(U, \hbox{\rm Bil}_{({\cal V},{\cal I})})$ 
 is given on sections (over open subsets of $U$) by $v\otimes v'\mapsto q(v+v')-q(v)-q(v').$
 Given a (not-necessarily symmetric!) bilinear form
 $b$, we also have the induced quadratic form $q_b$ given on sections by 
 $v\mapsto b(v\otimes v).$
 A global quadratic form may not be induced 
 from a global bilinear form, unless we assume something more, for e.g., 
 that the scheme is affine, or more 
 generally  that the sheaf cohomology group
 $\hbox{\rm H}^1(X, \hbox{\rm Alt}^2_{({\cal V}, {\cal I})})=0$.

\paragraph{The Generalised Clifford Algebra of Bichsel-Knus.}
Let$R$ be a commutative ring (with 1), $I$ an invertible $R$-module and
 $V$ a projective $R$-module. Let $L[I]:=R\oplus
\left(\bigoplus_{n> 0}(T^n(I) \oplus T^n(I^{-1}))\right)$ be the 
 Laurent-Rees algebra of  $I$, and define the $\mathbb Z$-gradation on the 
 tensor product of algebras $TV\otimes L[I]$ 
 by requiring elements of $V$ (resp. of $I$) 
 to be of degree one (resp. of degree two).
 Let $q:V\longrightarrow I$ be an $I$-valued quadratic form on $V.$
Following the definition of Bichsel-Knus \cite{Bichsel-Knus}, let
 $J(q, I)$ be the two-sided ideal of $TV\otimes L[I]$ generated by the set 
 $\{(x\otimes_{TV} x)\otimes 1_{L[I]}-1_{TV}\otimes q(x)\thinspace
|\thinspace x\in V\}$ and let the generalised Clifford algebra of $q$ be 
 defined by  $\widetilde{C}(V,q,I):=TV\otimes L[I]/J(q,I).$ This is an 
 $\mathbb Z$-graded algebra by definition. Let $C_n$ be the submodule of 
 elements of degree $n.$ Then $C_0$ is a subalgebra, playing the role of 
 the even Clifford algebra in the classical situation (i.e., $I=R$) and 
 $C_1$ is a $C_0$-bimodule. Bichsel and Knus baptize $C_0$ and $C_1$ 
 respectively as the {\em even Clifford algebra} and the {\em Clifford
  module} associated to the triple $(V,q,I).$
  The generalised Clifford algebra satisfies 
 an appropriate universal property which 
 ensures it behaves well functorially. Since $V$ is projective,  
 the canonical maps $V\longrightarrow \widetilde{C}(V,q,I)$ and 
 $L[I]\longrightarrow \widetilde{C}(V,q,I)$ are injective. For proofs of 
 these facts, see Sec.3 of \cite{Bichsel-Knus}. 
 If $(V,q,I)$ is an $I$-valued quadratic form on the vector bundle $V$
 over a scheme $X$, with $I$ a line bundle, then the above construction 
 may be carried out to define the generalised Clifford algebra bundle 
 $\widetilde{C}(V,q,I)$ which is an $\mathbb Z$-graded algebra bundle on 
 $X.$ 

\paragraph{Bourbaki's Tensor Operations with Values in a Line Bundle.}
 Let $R$ be a commutative ring and $L[I]:=R\oplus
\left(\bigoplus_{n> 0}(T^n(I) \oplus T^n(I^{-1}))\right)$ as above. 
 We denote by $\otimes_T$ (resp. by $\otimes_L$)
 the tensor product and by $1_T$ (resp. $1_L$) 
 the unit element in the algebra $TV$ (resp. in $L[I]$).
\begin{theorem}[{\rm with the above notations}]\label{Bourbaki}\   
\begin{description}
\item[(1)]
Let $q:V\longrightarrow I$ be an $I$-valued quadratic form on $V$ and 
$f\in \hbox{\rm Hom}_R(V, I).$ Then there exists an $R$-linear
 endomorphism $t_f$ of 
 the algebra $TV\otimes L[I]$ which is unique with respect to the first three
 of the following properties it satisfies:
\begin{description}
\item[(a)] $t_f(1_T\otimes \lambda)
=0\thinspace\forall\thinspace\lambda\in L[I];$
 \item[(b)]$t_f((x\otimes_T y)\otimes \lambda)
=y\otimes(f(x)\otimes_L\lambda)-(x\otimes_T1_T)t_f(y\otimes\lambda)$ 
for any $x\in V$, $y\in TV$, and $\lambda\in L[I];$
\item[(c)]following the definition of Bichsel-Knus \cite{Bichsel-Knus}, let
 $J(q, I)$ be the two-sided ideal of $TV\otimes L[I]$ generated by the set 
 $\{(x\otimes_T x)\otimes 1_L-1_T\otimes q(x)\thinspace
|\thinspace x\in V\}$; then $t_f(J(q, I))\subset J(q, I);$
\item[(d)]$t_f$ is homogeneous of degree $+1$ (except for elements which 
 it does not annihilate);
\item[(e)]recall from the definition of Bichsel-Knus \cite{Bichsel-Knus} 
that  the generalised Clifford algebra of $q$ is 
$\widetilde{C}(V,q,I):=TV\otimes L[I]/J(q,I)$;
 by (c) above, $t_f$ 
induces a $\mathbb Z$-graded endomorphism of degree $+1$ denoted by 
 $d_f^q:\widetilde{C}(V,q, I)\longrightarrow \widetilde{C}(V,q, I);$
\item[(f)]$t_f\circ t_f=0;$
\item[(g)]if $g\in \hbox{\rm Hom}_R(V,I)$, 
then $t_f\circ t_g + t_g\circ t_f=0;$
\item[(h)]if $\alpha\in \hbox{\rm End}_R(V)$, then 
 $t_f\circ (T(\alpha)\otimes \hbox{\rm Id}_{L[I]})=(T(\alpha)\otimes\hbox{\rm Id}_{L[I]})\circ t_{\alpha^*f}$ where $\alpha^*f\in \hbox{\rm Hom}_R(V, I)$ is
 defined by  $x\mapsto f(\alpha(x));$
\item[(i)]$t_f\equiv 0$ on the $R$-subalgebra of $TV\otimes L[I]$ 
generated by $\hbox{\rm kernel}(f)\otimes L[I].$ In fact, atleast when  
 $V$ is a projective $R$-module, the smallest $R$-subalgebra of
 $TV\otimes L[I]$ containing $\hbox{\rm kernel}(f)\otimes R.1_L$ is 
 $\hbox{\rm Image}(T(\hbox{\rm kernel}(f))\otimes R.1_L$ and $t_f$ vanishes 
 on this $R$-subalgebra. 
\end{description} 
\item[(2)] Let $q, q':V\longrightarrow I$ be two $I$-valued 
 quadratic forms whose
 difference is the quadratic form $q_b$ induced by 
 an $I$-valued 
  bilinear form $b\in\hbox{\rm Bil}_R(V, I)
:=\hbox{\rm Hom}_R(V\otimes_R V, I)$ 
i.e., $q'(x)-q(x)=q_b(x):=b(x,x)\thinspace \forall x\in V.$
 Further, for any $x\in V$ denote by $b_x$ the element of 
 $\hbox{\rm Hom}_R(V, I)$  given by 
 $y\mapsto b(x,y).$ Then there exists an $R$-linear automorphism 
 $\Psi_b$ of 
 $TV\otimes L[I]$
 which is unique with respect to the first three of the following 
 properties it satisfies: 
\begin{description}
\item[(a)]$\Psi_b(1_T\otimes \lambda)=(1_T\otimes\lambda)\thinspace \forall 
\lambda\thinspace \in L[I];$
\item[(b)]$\Psi_b((x\otimes_T y)\otimes\lambda)
=(x\otimes 1_L).\Psi_b(y\otimes\lambda)+t_{b_x}(\Psi_b(y\otimes\lambda))$ 
 for any $x\in V$, $y\in TV$ and $\lambda\in L[I];$ 
\item[(c)]$\Psi_b(J(q', I))\subset J(q, I);$ 
\item[(d)]by the previous property, $\Psi_b$ induces an
 isomorphism of $\mathbb Z$-graded 
 $R$-modules $$\psi_b:\widetilde{C}(V,q', I)\cong \widetilde{C}(V,q, I);$$
 in particular, given a 
 quadratic form $q_1:V\longrightarrow I$, since there always 
 exists an $I$-valued  bilinear form $b_1$ that induces $q_1$ (i.e.,
 such that $q_1=q_{b_1}$), setting $q'=q_1$, $q=0$ and $b=b_1$ in the above 
 gives an   
 $\mathbb Z$-graded linear
 isomorphism $\psi_{b_1}:\widetilde{C}(V,q_1, I)\cong 
 \widetilde{C}(V,0, I)=\Lambda(V)\otimes L[I];$ 
\item[(e)]$\Psi_b(T^{2n}V\otimes L[I])\subset
 \oplus_{(i\leq n)}(T^{2i}V\otimes L[I])$, 
 $\Psi_b(T^{2n+1}V\otimes L[I])\subset
 \oplus_{(\hbox{\rm\tiny odd }i\leq 2n+1)}(T^{i}V\otimes L[I])$, 
$\Psi_b(T^{2n}V\otimes I^{-n})\subset
 \oplus_{(i\leq n)}(T^{2i}V\otimes I^{-i})$ and 
 $\Psi_b(T^{2n+1}V\otimes I^{-n})\subset
 \oplus_{(\hbox{\rm\tiny odd }i\leq 2n+1)}(T^{i}V\otimes I^{\frac{1-i}{2}});$
\item[(f)]in particular, for $x,x'\in V$, $\Psi_b((x\otimes_T x')\otimes 1_L)
 =(x\otimes_T x')\otimes 1_L + 1_T\otimes b(x,x')$ 
 so that for $\psi_b:C_0(V,q_b, I)\cong C_0(V,0, I)
=\oplus_{n\geq 0}(\Lambda^{2n}(V)\otimes I^{-n})$ we have 
 $\psi_b(((x\otimes_T x')\otimes \zeta)\textrm{ mod }J(q_b, I))=(x\wedge x')
\otimes\zeta+ \zeta(b(x,x')).1$ for any $x,x'\in V$ and 
$\zeta\in I^{-1}\equiv\hbox{\rm Hom}_R(I, R);$  
\item[(g)]if $f\in \hbox{\rm Hom}_R(V, I)$, and $t_f$ is given by (1) above,
 then $\Psi_b\circ t_f=t_f\circ \Psi_b;$
\item[(h)]for $I$-valued bilinear forms $b_i$ on $V$,
 $\Psi_{b_1+b_2}=\Psi_{b_1}\circ 
 \Psi_{b_2}$ and $\Psi_0=\hbox{\rm Identity on } {TV\otimes L[I]};$ 
\item[(i)]for any $\alpha\in\hbox{\rm End}_R(V)$, $\Psi_b\circ
 (T(\alpha)\otimes \hbox{\rm Id}_{L[I]})=(T(\alpha)\otimes\hbox{\rm Id}_{L[I]})
\circ \Psi_{(b.\alpha)}$ where $(b.\alpha)(x,x'):=b(\alpha(x),\alpha(x'))\thinspace\forall x,x'\in V;$ 
\item[(j)]by property (h), 
one has a homomorphism of groups $(\hbox{\rm Bil}_R(V, I), +)
\longrightarrow (\hbox{\rm Aut}_R(TV\otimes L[I]), \circ):b\mapsto \Psi_b;$
 the associative
 unital monoid $(\hbox{\rm End}_R(V), \circ)$ acts on $\hbox{\rm Bil}_R(V, I)$ 
 on the right by $b'\leadsto b'.\alpha$ and acts on the 
 left (resp. on the right) of $\hbox{\rm End}_R(TV\otimes L[I])$ by
 $\alpha.\Phi:=(T(\alpha)\otimes\hbox{\rm Id}_{L[I]})\circ\Phi$
 (resp. by $\Phi.\alpha:=\Phi\circ (T(\alpha)\otimes\hbox{\rm Id}_{L[I]})$,
 and the  homomorphism
 $b\mapsto \Psi_b$ satisfies $\alpha.\Psi_{(b.\alpha)}=\Psi_b.\alpha;$ 
 the group $\hbox{\rm Aut}_R(V)=\hbox{\rm GL}_R(V)$ acts on
  the left of $\hbox{\rm Bil}_R(V, I)$ by 
 $g.b: (x,x')\mapsto b(g^{-1}(x), g^{-1}(x'))$ and on the left of 
 $\hbox{\rm Aut}_R(TV\otimes L[I])$ by conjugation via the natural 
 group homomorphism $\hbox{\rm GL}_R(V)\longrightarrow
 \hbox{\rm Aut}_R(TV\otimes L[I]): 
 g\mapsto T(g)\otimes\hbox{\rm Id}_{L[I]}$
 i.e., $g.\Phi:=(T(g)\otimes\hbox{\rm Id}_{L[I]})\circ \Phi\circ 
 (T(g^{-1})\otimes\hbox{\rm Id}_{L[I]})$, and the 
 homomorphism $b\mapsto \Psi_b$ is $\hbox{\rm GL}_R(V)$-equivariant: 
 $\Psi_{g.b}=g.\Psi_b.$ 
\end{description} 
\item[(3)] For a commutative $R$-algebra $S$ (with 1), let 
 $(q\otimes_R S), (q'\otimes_R S):
 (V\otimes_R S=:V_S)\longrightarrow (I\otimes_R S=:I_S)$ 
 be the $I_S$-valued quadratic forms induced from 
 the quadratic forms $q,q'$ of (2) above
 and $(b\otimes_R S)\in\hbox{\rm Bil}_S(V_S, I_S)$
 the $I_S$-valued $S$-bilinear form induced 
 from the bilinear form $b$ of (2) above. 
 Then as a result  of the uniqueness properties (2a)--(2c) satisfied 
 by $\Psi_b$ and $\Psi_{(b\otimes_R S)}$, the $S$-linear automorphisms
 $(\Psi_b\otimes_R S)$ and $\Psi_{(b\otimes_R S)}$ 
may be canonically  identified. In particular,
 the $\mathbb Z$-graded $S$-linear isomorphism 
 $(\psi_{b}\otimes_R S):\widetilde{C}(V_S,(q'\otimes_R\thinspace S), I_S)\cong\widetilde{C}(V_S,(q\otimes_R S), I_S)$ 
 induced from $\psi_b$ of
 (2d) above may be canonically identified with $\psi_{(b\otimes_R S)}.$ 
\end{description}
\end{theorem}

\begin{remark}\rm 
As mentioned in \cite{Bichsel-Knus}, F. van Oystaeyen has observed that 
 $L[I]$ is a faithfully-flat splitting for $I$, and the generalised 
 Clifford algebra $\widetilde{C}(V,q,I)$ is nothing but the ``classical'' 
 Clifford algebra of the triple $(V\otimes_R L[I], q\otimes_R\thinspace L[I], 
 I\otimes_R L[I]\equiv L[I])$ over $L[I].$  In the same vein, $I$-valued forms 
 (both multilinear and quadratic) on an $R$-module $V$ can be treated as 
 the usual ($L[I]$-valued) forms on $V\otimes_R L[I].$ 
 With this in mind, the proof of the above proposition follows from the 
 usual Bourbaki tensor operations with respect to $V\otimes_R L[I]$ on 
 $L[I].$  (See \S 9, Chap.9, \cite{Bourbaki} or 
 para.1.7, Chap.IV, \cite{Knus}, or Theorem 2.1 \cite{tevb-paper2} 
 for the ``classical'' Bourbaki operations). 
 However one needs to remember that the $\mathbb Z$-gradation on 
 $TV\otimes_R L[I]$ as defined above is different from the usual 
 ${\mathbb Z}_{\geq 0}$-gradation on $T(V\otimes_R L[I]).$  
\end{remark}


\paragraph{Tensoring by Symmetric Bilinear Bundles and
 Twisted Discriminant Bundles.}
Let $V,M$ be vector bundles on a scheme $X$ and let $I,J$ be line bundles on 
 $X.$ Let $q$ be a quadratic form on $V$ with values in $I$ and let 
 $b$ be a symmetric bilinear form on $M$ with values in $J.$
\begin{proposition}[{\rm with the above notations}]
\label{tensoring-with-sym-bil-module}\
\begin{description}
\item[(1)] The tensor product of $\vqi$ with $(M,b, J)$ gives a unique 
 quadratic bundle $(V\otimes M, 
 q\otimes b, I\otimes J).$ The quadratic form on $V\otimes M$ is given on 
 sections by $v\otimes m\mapsto q(v)\otimes b(m\otimes m)$
 and has associated bilinear 
 form $b_{q\otimes b}=b_q\otimes b.$
\item[(2)] When $M$ is a line bundle, $(M,b,J)$ is regular (=nonsingular)
  iff $(M, q_b, J)$  is 
 semiregular iff $b:M\otimes M\cong J$ is an 
 isomorphism (i.e., in other words  iff $(M,b,J)$ is a twisted 
 discriminant bundle). 
\item[(3)] Let $V$ be of odd rank and $(M, b, J)$ a twisted 
 discriminant bundle.
 Then $\vqi$ is semiregular iff
 $\vqi\otimes(M,b,J)=(V\otimes M, q\otimes b, I\otimes J)$ is semiregular. 
\end{description}
\end{proposition} 

\begin{proposition}\label{tensoring-with-disc-module}
Let $V$ and $V'$ be vector bundles of the same rank on the scheme $X$, 
 $(L,h, J)$ a twisted discriminant line bundle on $X$ and 
 $\alpha:V'\cong V\otimes L$ an isomorphism of bundles. 
\begin{description}
\item[(1)] Over any open subset
 $U\hookrightarrow X$, given a bilinear form 
 $b'\in\Gamma\left(U, \hbox{\rm Bil}_{(V',I)}\right)$,
 we can define a 
 bilinear form $b\in
\Gamma\left(U, \hbox{\rm Bil}_{(V, I\otimes J^{-1})} \right)$ 
 using $\alpha$ and $h$ as follows: we let 
$b:=(b'\otimes J^{-1})\circ{(\zeta_{(\alpha, h)})}^{-1}$
 where $\zeta_{(\alpha,h)}:V'\otimes V'\otimes J^{-1}\cong
 V\otimes V$ is the linear isomorphism given
 by the composition of the following natural morphisms:
\begin{center}\begin{tabular}{c}
$V'\otimes V'\otimes J^{-1}
\stackrel{\alpha\otimes\alpha\otimes\hbox{\rm\tiny Id}
\thinspace(\cong)}{\longrightarrow}
V\otimes L\otimes V\otimes L\otimes J^{-1}
\stackrel{\hbox{\rm\tiny SWAP}(2,3)\thinspace(\equiv)}
{\longrightarrow}
V\otimes V\otimes L^2\otimes J^{-1}
\stackrel{\hbox{\rm\tiny Id}\otimes h\otimes\hbox{\rm\tiny Id}
\thinspace(\cong)}{\longrightarrow}$\\
$\stackrel{\hbox{\rm\tiny Id}\otimes h\otimes\hbox{\rm\tiny Id}
\thinspace(\cong)}{\longrightarrow}V\otimes V\otimes J\otimes J^{-1}
\stackrel{\hbox{\rm\tiny CANON}\thinspace(\equiv)}{\longrightarrow}
V\otimes V.$ \end{tabular}\end{center}
Then the association $b'\mapsto b$ induces linear isomorphisms 
 shown by vertical upward arrows in the following diagram of associated 
 locally-free sheaves (with 
 exact rows) making it commutative: 
$$\begin{CD}
0 @>>> \hbox{\rm Alt}^2_{({\cal V},{\cal I}\otimes {\cal J}^{-1})}
 @>>> \hbox{\rm Bil}_{({\cal V}, {\cal I}\otimes {\cal J}^{-1})}
 @>>> \hbox{\rm Quad}_{({\cal V}, {\cal I}\otimes {\cal J}^{-1})} @>>> 0\\
 & & @A{\cong}AA  @A{\cong}AA  @A{\cong}AA & & \\
0 @>>> \hbox{\rm Alt}^2_{({\cal V}', {\cal I})} @>>> \hbox{\rm Bil}_{({\cal V}',{\cal I})} @>>>
 \hbox{\rm Quad}_{({\cal V}', {\cal I})} @>>> 0.
\end{CD}$$ 
Therefore one also has the following commutative diagram of vector bundle 
 morphisms with the vertical upward arrows being isomorphisms: 
$$\begin{CD}
 \hbox{\rm Alt}^2_{(V, I\otimes J^{-1})} 
@>{\hbox{\rm closed}}>{\hbox{\rm immersion}}> 
\hbox{\rm Bil}_{(V, I\otimes J^{-1})}
 @>{\hbox{\rm locally}}>{\hbox{\rm trivial}}>
 \hbox{\rm Quad}_{(V, I\otimes J^{-1})} \\
 @A{\cong}AA  @A{\cong}AA @A{\cong}AA \\
\hbox{\rm Alt}^2_{(V', I)} 
@>{\hbox{\rm closed}}>{\hbox{\rm immersion }}>
 \hbox{\rm Bil}_{(V', I)}
 @>{\hbox{\rm locally}}>{\hbox{\rm trivial}}>
 \hbox{\rm Quad}_{(V', I)} 
\end{CD}$$ 
\item[(2)]
 Let $b'\in
\Gamma\left(X, \hbox{\rm Bil}_{(V', I)}\right)$ be a global 
 bilinear form and let it induce 
$b\in\Gamma\left(X, \hbox{\rm Bil}_{(V, I\otimes J^{-1})}\right)$
 via $\alpha$ and $h$  as defined in (1) above.
 Let $\Psi_{b'}\in\hbox{\rm Aut}_{{\cal O}_X}(TV'\otimes L[I])$ 
(resp. $\Psi_{b}\in
\hbox{\rm Aut}_{{\cal  O}_X}(TV\otimes L[I\otimes J^{-1}])$) 
be the $\mathbb Z$-graded linear isomorphism   
 induced by 
 $b'$ (resp. by $b$) defined locally (and hence globally) 
 as in (2), Theorem \ref{Bourbaki} above.
 Let $Z_{(\alpha, h)}:\oplus_{n\geq 0}(T^{2n}_{{\cal  O}_X}(V')\otimes I^{-n})
\cong\oplus_{n\geq 0}(T^{2n}_{{\cal  O}_X}(V)\otimes I^{-n}\otimes J^n)$
 be the ${\cal  O}_X$-algebra isomorphism induced via the isomorphism
 $\zeta_{(\alpha,h)}$ defined in (1) above. Then, taking into account 
 (2e), Theorem \ref{Bourbaki},  
 the following diagram commutes:
$$\begin{CD}
\oplus_{n\geq 0}(T^{2n}_{{\cal  O}_X}(V')\otimes I^{-n}) 
@>{Z_{(\alpha, h)}}>{\cong}>
\oplus_{n\geq 0}(T^{2n}_{{\cal  O}_X}(V)\otimes I^{-n}\otimes J^n)\\
@V{\Psi_{b'}}V{\cong}V  @V{\cong}V{\Psi_b}V \\
\oplus_{n\geq 0}(T^{2n}_{{\cal  O}_X}(V')\otimes I^{-n}) 
@>{\cong}>{Z_{(\alpha, h)}}>
\oplus_{n\geq 0}(T^{2n}_{{\cal  O}_X}(V)\otimes I^{-n}\otimes J^n)
\end{CD}$$ 
thereby inducing by (2d), Theorem \ref{Bourbaki} 
the following commutative diagram of ${\cal  O}_X$-linear 
 isomorphisms
$$\begin{CD}
C_0(V', q_{b'}, I) @>{\hbox{\rm via }Z_{(\alpha, h)}}>{\cong}>C_0(V, q_b, I\otimes J^{-1})\\
@V{\psi_{b'}}V{\cong}V @V{\cong}V{\psi_b}V \\
\oplus_{n\geq 0}(\Lambda^{2n}_{{\cal  O}_X}(V')\otimes I^{-n})
 @>{\cong}>{\hbox{\rm via }Z_{(\alpha, h)}}>
\oplus_{n\geq 0}(\Lambda^{2n}_{{\cal  O}_X}(V)\otimes I^{-n}\otimes J^n)
\end{CD}$$ 
\item[(3)]
 Let $b$ and $b'$ be as in (2) above. Then $\alpha:V'\cong V\otimes L$ 
 induces an isometry 
 of bilinear form bundles $\alpha:(V',b', I)\cong(V,b, I\otimes J^{-1})
\otimes(L,h, J)$ and also an isometry of the induced 
 quadratic bundles $\alpha: (V',q_{b'}, I)\cong 
 (V, q_b, I\otimes J^{-1})\otimes (L,h, J).$
 Moreover, if we are just given a global $I\otimes J^{-1}$-valued quadratic 
 form $q$ on $V$ (resp. an $I$-valued $q'$ on $V'$),
 then we may define the global $I$-valued quadratic 
 form $q'$ on $V'$ (resp. $I\otimes J^{-1}$-valued $q$ on $V$)
 via $q':=(q\otimes h)\circ \alpha$
 (resp. via 
$q:=(q'\circ\alpha^{-1})\otimes{(h^\vee)}^{-1})$) 
 and  
 again $\alpha:(V',q', I)\cong (V,q, I\otimes J^{-1})\otimes(L,h, J)$
 becomes an
 isometry of quadratic 
 bundles.   
\end{description} 
\end{proposition}  
 
\begin{proposition}\label{simil-induces-iso-of-even-cliff}
Let $g:\vqi\thinspace{\cong}_{l}\thinspace\vprqpri$
 be an $I$-similarity with multiplier
 $l\in\Gamma(X, {\cal  O}_X^*).$ 
\begin{description}
\item[(1)] There exists a unique isomorphism of ${\cal  O}_X$-algebra bundles
$C_0(g,l, I):C_0(V,q, I)\cong C_0(V',q', I)$ such that  
 $C_0(g,l, I)(v.v'.s)=g(v).g(v').l^{-1}s$ on sections $v, v'$ of $V$ and 
 $s$ of $I^{-1}.$
\item[(2)] There exists a unique vector bundle isomorphism 
$C_1(g,l, I):C_1(V,q, I)\cong C_1(V', q', I)$ such that
\begin{description}
\item[(a)] 
 $C_1(g,l, I)(v.c)=g(v).C_0(g,l, I)(c)$ and 
\item[(b)] $C_1(g,l, I)(c.v)=C_0(g,l,I)(c).g(v)$ 
\end{description} 
for any section $v$ of $V$ and any section $c$ of $C_0(V,q).$
 Thus $C_1(g,l, I)$ is $C_0(g,l, I)$-semilinear.  
\item[(3)]
 If $g_1:\vprqpri\thinspace{\cong}_{l_1}\thinspace(V'',q'',I)$ is another 
  similarity with multiplier $l_1$, then the composition 
 $g_1\circ g:\vqi\thinspace{\cong}_{ll_1}\thinspace(V'',q'', I)$ is 
 also a similarity with multiplier given by the product of 
 the multipliers. Further $C_i(g_1\circ g,ll_1,I)=
 C_i(g_1,l_1, I)\circ C_i(g,l, I)$ for $i=0,1.$ 
\end{description}
\end{proposition}
\noindent 
 A local computation 
 shows that tensoring by a twisted 
 discriminant bundle amounts to (locally) applying 
 a similarity. In this case also one gets a global
  isomorphism of even Clifford 
 algebras: 

\begin{proposition}\label{isom-covq-covlqh}
Let $\vqi$ be a quadratic bundle on a scheme $X$ and $(L,h, J)$ be a twisted 
 discriminant bundle. There exists a unique isomorphism of algebra bundles 
 $$\gamma_{(L,h,J)}:C_0\left(\vqi\otimes(L,h,J)\right)\cong C_0(V,q,I)$$ 
 given by 
$\gamma_{(L,h,J)}\left((v\otimes\lambda).
(v'\otimes \lambda').(s\otimes t)\right)=
t(h(\lambda\otimes \lambda'))v.v'.s$ 
for any  sections $v,v'$ of $V$, $\lambda,\lambda'$ of $L$, 
 $s$ of $I^{-1}\equiv I^{\vee}$ and $t$ of $J^{-1}\equiv J^{\vee}.$   
\end{proposition}

\section{\label{sec3} Proof of Injectivity: Theorems \ref{lifting-of-isomorphisms} \& \ref{lifting-of-automorphisms}}

\paragraph{Proof of Prop.\ref{transfer-to-lambda2}:}
Start with an isomorphism of algebra-bundles $\phi:\covqi\cong\covprqpripr.$
Let ${\{U_i\}}_{i\in \mathcal{I}}$
 be an affine open covering of $X$ (which may 
also be chosen so as to trivialise some or any of 
 the involved bundles if needed). Choose bilinear forms 
 $b_i\in\Gamma(U_i, \hbox{\rm Bil}_{(V, I)})$ 
and $b'_i\in\Gamma(U_i, \hbox{\rm Bil}_{(V', I')})$ 
such that $q|U_i=q_{b_i}$ and $q'|U_i=q_{b'_i}$ for each $i\in \mathcal{I}.$ 
 By (2d), Theorem \ref{Bourbaki}, we have  isomorphisms of vector bundles 
$\psi_{b_i}$ and 
 $\psi_{b'_i}$, which preserve 1 by (2a) of the same Theorem, 
  and we define the isomorphism of vector bundles $\phi_{\Lambda_i^{ev}}$ 
 so as to make the following diagram commute: 
$$\begin{CD}
C_0(V,q, I)|U_i @>{\phi|U_i}>{\cong}> C_0(V', q', I')|U_i\\
@V{\psi_{b_i}}V{\cong}V @V{\cong}V{\psi_{b'_i}}V\\
({\cal O}_X\oplus\Lambda^2(V)\otimes I^{-1})|U_i 
@>{\cong}>{\phi_{\Lambda_i^{ev}}}> 
({\cal O}_X\oplus\Lambda^2(V')\otimes {(I')}^{-1})|U_i 
\end{CD}$$
The linear isomorphism $\phi_{\Lambda_i^{ev}}$ preserves 1  and 
 therefore it induces a linear isomorphism from 
$(\Lambda^2(V)\otimes I^{-1})|U_i$ to 
 $(\Lambda^2(V')\otimes {(I')}^{-1})|U_i$,
 which we denote by  ${({\phi_{\Lambda^2}})}_i.$ 
Observe that ${({\phi_{\Lambda^2}})}_i$ is independent of the choice of 
 the bilinear forms $b_i$ and $b'_i.$ For, replacing these respectively 
 by $\widehat{b_i}$ and $\widehat{b'_i}$, it follows from (2f),
 Theorem \ref{Bourbaki}, that $\psi_{b_i}\circ{(\psi_{\widehat{b_i}})}^{-1}$ 
 (resp. $\psi_{b'_i}\circ{(\psi_{\widehat{b'_i}})}^{-1}$) followed by the 
 canonical projection onto $(\Lambda^2(V)\otimes I^{-1})|U_i$ 
(resp. onto $(\Lambda^2(V')\otimes {(I')}^{-1})|U_i$) 
 is the same as the projection itself. By this observation, it is also 
 clear that the isomorphisms 
${\{{({\phi_{\Lambda^2}})}_i\}}_{i\in \mathcal{I}}$
 agree on (any open affine subscheme of, and hence on all of)
 any intersection $U_i\cap U_j.$ Therefore they glue to give a global 
 isomorphism of vector bundles 
$\phi_{\Lambda^2}:\Lambda^2(V)\otimes I^{-1}\cong
\Lambda^2(V')\otimes {(I')}^{-1}.$ 
{\bf Q.E.D, Prop.\ref{transfer-to-lambda2}.}

 \paragraph{Reduction of Proof of Injectivity of Theorem \ref{bijectivity} to 
 Theorem \ref{lifting-of-isomorphisms}.}
 We start with an isomorphism of algebra-bundles
 $\phi:\covqi\cong\covprqpripr$, construct the isomorphism of vector
 bundles $\phi_{\Lambda^2}:\Lambda^2(V)\otimes I^{-1}
\cong\Lambda^2(V')\otimes {(I')}^{-1}$
 and keep the notations
 introduced in the proof of Prop.\ref{transfer-to-lambda2}. 
 Firstly we deduce a linear isomorphism 
  $\hbox{\rm det}{({(\phi_{\Lambda^2})}^\vee)}^{-1}:
\hbox{\rm det}({(\Lambda^2(V)\otimes I^{-1})}^\vee)\cong
\hbox{\rm det}({(\Lambda^2(V')\otimes {(I')}^{-1})}^\vee).$ 
 Since $V$ and  $V'$ are of rank 3,  
  there are canonical isomorphisms
 $\eta:\Lambda^2(V)\equiv V^\vee\otimes\hbox{\rm det}(V)$
 and $\eta':\Lambda^2(V')\equiv (V')^\vee\otimes\hbox{\rm det}(V').$
It follows therefore that if we set $L:=\hbox{\rm det}(V')
\otimes{(\hbox{\rm det}(V))}^{-1}$ and $J:=I'\otimes I^{-1}$
 then we get a twisted discriminant line bundle $(L\otimes J^{-1},h, J)$
 and a vector 
 bundle isomorphism $\alpha:V'\cong V\otimes (L\otimes J^{-1}).$

 Now for each $i\in \mathcal{I}$, 
 the bilinear form 
 $b_i\in\Gamma(U_i, \hbox{\rm Bil}_{(V,I)})$
 induces, via $\alpha|U_i$ and $(LJ^{-1},h, J)|U_i$ and
 (1), Prop.\ref{tensoring-with-disc-module}, a bilinear form 
 $b''_i\in\Gamma(U_i, \hbox{\rm Bil}_{(V', IJ)}).$
 By (3) of the same Proposition, 
 over each $U_i$  we get an isometry  of bilinear form bundles 
$\alpha|U_i:(V'|U_i, b''_i, IJ|U_i)\cong (V|U_i, b_i, I|U_i)
\otimes (LJ^{-1}, h, J)|U_i$ 
 and also an isometry of quadratic bundles
$\alpha|U_i:(V'|U_i, q_{b''_i}, IJ|U_i)\cong 
(V|U_i, q_{b_i}=q|U_i, I|U_i)\otimes (LJ^{-1}, h, J)|U_i.$

 On the other hand,
 by an assertion in (3), Prop.\ref{tensoring-with-disc-module}, we could also 
 define the global
 quadratic bundle $(V',q'', IJ)$ using $(V,q, I)$, $\alpha$ and
 $(LJ^{-1},h, J)$, 
 so that we have an isometry of quadratic bundles 
$\alpha:(V',q'', IJ)\cong(V,q, I)\otimes(LJ^{-1},h, J).$
 It follows therefore that 
 the $q_{b''_i}$ glue to give $q''.$ 
Notice that in general 
 the $b''_i$ (resp. the $b_i$) need not glue to give a global bilinear 
 form $b''$ (resp. $b$) such that $q_{b''}=q''$ (resp. $q_b=q$). 
 By (1), Prop.\ref{simil-induces-iso-of-even-cliff}, there exists 
 a  unique isomorphism of algebra bundles
$$C_0(\alpha,1,IJ):C_0(V',q'', IJ)\cong C_0\left((V,q, I)
\otimes(LJ^{-1},h, J)\right)$$ and by 
 Prop.\ref{isom-covq-covlqh} 
 we have a unique isomorphism of algebra bundles 
 $$\gamma_{(LJ^{-1},h,J)}:C_0\left(\vqi\otimes(LJ^{-1},h,J)\right)
\cong\covqi.$$
 Therefore the composition 
 of the following sequence of isomorphisms of algebra bundles on X
\begin{center}\begin{tabular}{c}
$C_0(V',q'', I')\stackrel{\hbox{\rm\tiny using }I'\cong IJ}{\cong}
C_0(V',q'', IJ)\stackrel{C_0(\alpha,1,IJ)\thinspace{(\cong)}}{\longrightarrow}
C_0\left((V,q,I)\otimes(LJ^{-1},h,J)\right)\stackrel{\gamma_{(LJ^{-1},h,J)}
\thinspace{(\cong)}}{\longrightarrow}$\\
$\stackrel{\gamma_{(LJ^{-1},h,J)}
\thinspace{(\cong)}}{\longrightarrow}
\covqi\stackrel{\phi\thinspace{(\cong)}}{\longrightarrow}\covprqpripr$
\end{tabular}\end{center}
is an element of $\hbox{\rm Iso}[C_0(V',q'', I'),C_0(V',q',I')]$,
 which, granting 
 Theorem \ref{lifting-of-isomorphisms}, is induced by a similarity in  
 $(g,l)\in\hbox{\rm Sim}[(V',q'', I'),(V',q',I')].$ Hence we would have 
 $$g:(V',q'',I')\cong(V',q',I')\otimes({\cal O}_X,\textrm{ mult.by }l^{-1}, 
 {\cal O}_X)$$ where $l\in\Gamma(X, {\cal O}_X^*).$
 This combined with the fact that 
   $(V,q,I)$ and $(V',q'',I')\cong(V',q'',IJ)$ are isomorphic 
 upto the twisted discriminant bundle $(LJ^{-1},h, J)$
 by the construction above would imply 
 that $(V,q,I)$ and $(V',q',I')$ also differ by a
twisted  discriminant 
 bundle. Therefore the proof of the injectivity asserted in  
 Theorem \ref{bijectivity} reduces 
 to the proof of 
 Theorem \ref{lifting-of-isomorphisms}. 

\paragraph{Reduction of Theorem \ref{lifting-of-isomorphisms}
 to the Case when $I$ is free:}
\label{reduction-to-free-case}
 For a similarity $g$ with multiplier $l$, we have $C_0(g,l,I)$
 given by (1), Prop.\ref{simil-induces-iso-of-even-cliff}, so that we 
 may define the map $$\simvqvqpri \longrightarrow \isomcovqcovqpri:
g\mapsto C_0(g,l,I).$$   
 The equality $\hbox{\rm det}(C_0(g,l, I))
=\hbox{\rm det}\left((C_0(g,l, I))_{\Lambda^2}\right)
=l^{-3}\hbox{\rm det}^2(g)$ was shown to hold (locally, hence globally) 
 (1), Lemma 3.11, \cite{tevb-paper2}.
 Thus $\isomvqvqpri$ and $\sisomvqvqpri$
 are respectively mapped 
 into  
 $\isomcovqcovqprprimei$ and $\sisomcovqcovqpri$ as claimed.
  We start with an isomorphism of algebra-bundles $\phi:\covqi\cong\covqpri$,
 which by Prop.\ref{transfer-to-lambda2} leads to the automorphism of vector
 bundles $\phi_{\Lambda^2}\in\hbox{\rm Aut}(\Lambda^2(V)\otimes I^{-1}).$
Firstly, define the global bundle automorphism 
 $g'\in\hbox{\rm GL}\left(V\otimes{(\hbox{\rm det}(V))}^{-1}\otimes I\right)$
 so that the following diagram commutes:
$$\begin{CD}
{(\Lambda^2(V))}^\vee\otimes I 
@>{{({({\phi_{\Lambda^2}})}^\vee)}^{-1}}>{\cong}> 
{(\Lambda^2(V))}^\vee\otimes I\\
@V{{(\eta^\vee)}^{-1}\otimes I}V{\equiv}V 
@V{\equiv}V{{(\eta^\vee)}^{-1}\otimes I}V\\
V\otimes{(\hbox{\rm det}(V))}^{-1}\otimes I
 @>{\cong}>g'> V\otimes{(\hbox{\rm det}(V))}^{-1}\otimes I
\end{CD}$$ where 
 $\eta:\Lambda^2(V)\equiv V^\vee\otimes\hbox{\rm det}(V)$
 is the canonical isomorphism (since $V$ is of rank 3). Now let 
 $g\in\hbox{\rm GL}(V)\stackrel{\cong}{\longleftarrow}\hbox{\rm GL}
(V\otimes{(\hbox{\rm det}(V))}^{-1}\otimes I)$ be the image of $g'$ i.e., 
 the image of $g'\otimes \hbox{\rm det}(V)\otimes I^{-1}$
 under the canonical identification 
 $\hbox{\rm GL}(V\otimes{(\hbox{\rm det}(V))}^{-1}\otimes I\otimes\hbox{\rm det}(V)\otimes I^{-1})
\equiv \hbox{\rm GL}(V).$ 
 Next, let $l\in\Gamma(X, {\cal O}_X^*)$ be a global section such that 
 $\gamma(l):=(l^3).\hbox{\rm det}(\phi_{\Lambda^2})$ has a square root in 
 $\Gamma(X, {\cal O}_X^*).$ For example, we have the following special 
 cases when this is true: 
\begin{description}
\item[Case 1.]  \label{cases}
 If $\hbox{\rm det}(\phi_{\Lambda^2})$ is itself a square, 
set $l:=1.$ If further $\hbox{\rm det}(\phi_{\Lambda^2})=1$, set 
 $\sqrt{\gamma(l)}=1$, otherwise let $\sqrt{\gamma(l)}$ denote any 
 fixed square root of $\hbox{\rm det}(\phi_{\Lambda^2}).$
\item[Case 2.] If $\hbox{\rm det}(\phi_{\Lambda^2})$ is not 
 a square, given an integer $k$,
 take $l={(\hbox{\rm det}(\phi_{\Lambda^2}))}^{2k+1}$ and
 let $\sqrt{\gamma(l)}$ denote any fixed square root of 
 ${(\hbox{\rm det}(\phi_{\Lambda^2}))}^{6k+4}.$
\end{description}
 For each integer 
 $k$, we now associate to $\phi$ the element \label{glphi-def}
 $g_l^\phi:=(l^{-1}\sqrt{\gamma(l)}\thinspace)g$ with 
 $g$ as defined above. Suppose we  
 show the following locally for the Zariski toplogy on $X$  
 (more precisely, for each open subscheme of $X$ over which 
 $V$ and $I$ are free): 
\begin{description}
\item[(1)]  
 that $g_l^\phi$ is an $I$-similarity from $(V,q,I)$ to $(V,q',I)$
 with multiplier $l$; 
\item[(2)] that $g_l^\phi$  
 induces $\phi$ i.e., with the notations of (1), 
Prop.\ref{simil-induces-iso-of-even-cliff}, 
 that $C_0(g_l^\phi, l, I)=\phi$; 
\item[(3)] that $\hbox{\rm det}(g_l^\phi)=\sqrt{\gamma(l)}$
 so that 
 $\hbox{\rm det}^2(g_l^\phi)
=\hbox{\rm det}(\phi_{\Lambda^2})$ 
when $\hbox{\rm det}(\phi_{\Lambda^2})$ is itself a square
 and 
\item[(4)]
 that the  map 
$\sisomvqvqpri\longrightarrow\sisomcovqcovqpri$
 is injective.
\end{description}
 It would follow then that these statements are also true globally. The 
 maps $s_{2k+1}:\phi\mapsto g_l^\phi$ with $l$ as in Case 2 and 
 $s': \phi\mapsto g_l^\phi$ with $l$ as in Case 1 will then give the 
 sections to the maps (which would imply their surjectivities) 
as mentioned in Theorem \ref{lifting-of-isomorphisms}.
 But these maps 
  are  not necessarily multiplicative since a computation reveals that 
 if $\phi_i\in\hbox{\rm Iso}[C_0(V,q_i,I),C_0(V,q_{i+1},I)]$ 
is associated to $g_{l_i}^{\phi_i}\in\hbox{\rm Sim}[(V,q_i,I),(V,q_{i+1},I)]$, 
 and $\phi_2\circ\phi_1$ to $g_{l_{21}}^{\phi_2\circ\phi_1}$, then 
 $g_{l_{21}}^{\phi_2\circ\phi_1}=\delta g_{l_2}^{\phi_2}\circ 
g_{l_1}^{\phi_1}$ for 
 $\delta\in\mu_2(\Gamma(X, {\cal O}_X))$ because of the ambiguity in 
 the initial global choices of square roots for 
 $\gamma(l_i)$ and $\gamma(l_{21}).$  However 
 this can be remedied as follows. For any given 
 $\phi\in\hbox{\rm Iso}[C_0(V,q,I),C_0(V,q',I)]$,
 irrespective of whether or not
 $\hbox{\rm det}(\phi_{\Lambda^2})$ is a square,
 take 
$$l={(\hbox{\rm det}(\phi_{\Lambda^2}))}^{2k+1}, 
 \gamma(l)=l^3\hbox{\rm det}(\phi_{\Lambda^2}),
  \sqrt{\gamma(l)}:={(\hbox{\rm det}(\phi_{\Lambda^2}))}^{3k+2}\hbox{ and }
 s_{2k+1}^{+}(\phi):=g_l^\phi=\left(l^{-1}\sqrt{\gamma(l)}\right)g.$$
 Then it is clear that each $s_{2k+1}^{+}$ is multiplicative with the 
 properties as claimed in the statement.
 We thus reduce the proof of Theorem \ref{lifting-of-isomorphisms} 
 to the case when the rank 3 vector bundle 
 $V$ and the line bundle $I$ are free. 
 More generally,  even for the case $V$ not necessarily free but 
 $I={\cal O}_X$, this was 
 treated in Theorem 1.3 of \cite{tevb-paper2}.
{\bf Q.E.D., Theorem \ref{lifting-of-isomorphisms} \& injectivity of 
 Theorem \ref{bijectivity}.}

\paragraph{Proof of Theorem \ref{lifting-of-automorphisms}.}
 (The case $I={\cal O}_X$ was treated in Theorem 1.4 of \cite{tevb-paper2}).
 Taking 
 $q'=q$ in Theorem \ref{lifting-of-isomorphisms} gives the commutative 
 diagram of groups and homomorphisms as asserted in the statement of 
 the theorem. We continue with the notations above. 
 For $g\in\govqi$ with multiplier $l$, that 
 the equality $\hbox{\rm det}(C_0(g,l, I))
=\hbox{\rm det}\left((C_0(g,l, I))_{\Lambda^2}\right)
=l^{-3}\hbox{\rm det}^2(g)$ holds (locally, hence globally) was 
 shown in (1),
 Lemma 3.11, \cite{tevb-paper2}.
 Assertion (2) of the same Lemma shows the following (locally, and 
 hence globally):  
 if $C_0(g,l, I)$ is the identity on 
$\covqi$, then 
$g=l^{-1}\hbox{\rm det}(g).\hbox{\rm Id}_V$,
 and 
 further if $g\in\ovqi$ then 
$g=\hbox{\rm det}(g).\hbox{\rm Id}_V$ with $\hbox{\rm det}^2(g)=1.$ 
 This gives 
 exactness at $\govqi$ and at $\ovqi.$ 

We proceed to prove assertion (b). 
 Let $\phi\in\autcovqi$, and consider the self-similarity
 $s_{2k+1}^{+}(\phi)=g_l^\phi$ with multiplier
 $l=\hbox{\rm det}(\phi)^{2k+1}.$ 
 For the moment, assume that $V$ and $I$ are trivial over 
 $X.$ Fix a
  global basis $\{e_1,e_2,e_3\}$ for $V$ 
  and set $e'_i=g_l^\phi(e_i).$ 
  It follows from Kneser's  definition 
 of the half-discriminant $d_0$---see formula (3.1.4), Chap.IV, 
\cite{Knus}---that $d_0(q, \{e_i\})
=d_0(q, \{e'_i\})\thinspace\hbox{\rm det}^2(g_l^\phi).$ 
 Since we have 
 $g_l^\phi.q=l^{-1}q$, a simple computation shows that $d_0(q, \{e'_i\})=
l^{3}d_0(q,\{e_i\}).$ The hypothesis that $q\otimes\kappa(x)$ is 
 semiregular means that the image of the 
 element $d_0(q,\{e_i\})\in\Gamma(X,{\cal O}_X)$ in $\kappa(x)$ is nonzero. 
 Since $X$ is integral, we therefore deduce that 
$\hbox{\rm det}^2(g_l^\phi)=l^{-3}.$ On the other hand, we know that 
 $\hbox{\rm det}^2(g_l^\phi)l^{-3}=\hbox{\rm det}(\phi).$ It follows 
 that $\hbox{\rm det}^{12k+7}(\phi)
=1\thinspace\forall\thinspace k\in\mathbb Z$,
 which forces $\hbox{\rm det}(\phi)=1.$ In general, even if $V$ and 
 $I$ are not necessarily trivial, since this equality holds 
 over a covering of $X$ which trivialises both $V$ and $I$, it 
 also holds over all of $X.$
{\bf Q.E.D., Theorem \ref{lifting-of-automorphisms}.}

\section{\label{sec4} Proof of Surjectivity: Theorems \ref{structure-of-specialisation} \& \ref{bilinear-forms-as-specialisations}}

\paragraph{\bf Proof of Theorem \ref{bilinear-forms-as-specialisations}.}
We adopt the notations introduced just before 
  Theorem \ref{bilinear-forms-as-specialisations}.
 Let $T$ be an $X$-scheme. Given a bilinear form $b\in\hbox{\rm Bil}_{(V,I)}(T)
=\Gamma(T, \hbox{\rm Bil}_{(V_T, I_T)})$,
 consider the linear isomorphism $\psi_{b}:C_0(V_T, q_b, I_T)
\cong {\cal O}_T\oplus\Lambda^2(V_T)\otimes {(I_T)}^{-1}=W_T$ of 
 (2d), Theorem \ref{Bourbaki}. Let $A_b$ denote the algebra bundle structure
 on $W_T$ with unit $w_T=1$ induced via $\psi_b$ 
 from the even Clifford algebra $C_0(V_T, q_b, I_T)$. By definition, 
 $A_b\in\aaWw(T)$ and we get a map of $T$-valued points 
 $$\Upsilon'(T):\hbox{\rm Bil}_{(V,I)}(T)
\longrightarrow \aaWw(T): b\mapsto A_b.$$
 This is functorial in $T$ because of (3), Theorem \ref{Bourbaki}, and 
 hence defines an $X$-morphism $\Upsilon':\hbox{\rm Bil}_{(V,I)}
\longrightarrow 
 \aaWw.$ The morphism $\Upsilon'$ is $\hbox{\rm GL}_V$-equivariant 
 due to (2j), Theorem \ref{Bourbaki}. Notice that the schemes 
 $\hbox{\rm Bil}_{(V,I)}, \hbox{\rm Bil}^{sr}_{(V,I)}$ and $\aaWw$ are 
 well-behaved relative to $X$ with respect to base-change. In fact, 
 so are $\azuWw$ and $\spazuWw$, as may be recalled from Theorems 
 3.4 and 3.8, Part A, \cite{tevb-paper1}. 
In view of these observations, by taking a trivialisation for $I$ over 
 $X$, we may reduce to the case when $I$ is trivial. But this case was 
 treated in Theorem 1.9 of \cite{tevb-paper2}. 
{\bf Q.E.D., Theorem \ref{bilinear-forms-as-specialisations}.}

\paragraph{Proofs of assertions in
  (a), Theorem \ref{structure-of-specialisation}
 and Surjectivity part of Theorem \ref{bijectivity}.} 
Let $W$ be the rank 4 vector bundle underlying 
 the specialised algebra $A$ and $w\in\Gamma(X, W)$
 be the global section corresponding to 
 $1_A.$ We choose an affine open covering 
${\{U_i\}}_{i\in \mathcal{I}}$ of $X$ such that $W|U_i$ is trivial and 
 $w|U_i$ is part of a global basis $\forall\thinspace i.$  
 Therefore on the one hand, 
for each $i\in \mathcal{I}$, we can find a linear isomorphism
 $\zeta_i:\Lambda^{even}\left({\cal O}_X^{\oplus 3}|U_i\right)\cong 
 W|U_i$ taking $1_{\Lambda^{even}}$ onto $w|U_i.$ The $(w|U_i)$-unital
 algebra structure  $A|U_i$ induces via $\zeta_i$ an algebra structure $A_i$ on
 $\Lambda^{even}\left({\cal O}_X^{\oplus 3}|U_i\right)$ (so that $\zeta_i$
 becomes an algebra isomorphism). Recall that $A_i$ is also 
 a specialised algebra structure by Theorem 3.8, Part A, \cite{tevb-paper1}.  
 Hence by  Theorem \ref{bilinear-forms-as-specialisations} 
 applied to $X=U_i$, $V={\cal O}_{U_i}^{\oplus 3}$ and $I={\cal O}_{U_i}$, 
  we can also
 find an ${\cal O}_{U_i}$-valued quadratic form 
 $q_i$ on ${\cal O}_X^{\oplus 3}|U_i$ induced from a bilinear form 
 $b_i$ so that the algebra structure $A_i$ is precisely the one 
 induced by the linear isomorphism 
 $\psi_{b_i}:C_0\left({\cal O}_X^{\oplus 3}|U_i, q_i\right)
\cong\Lambda^{even}\left({\cal O}_X^{\oplus 3}|U_i\right)$ 
given by (2d) of Theorem \ref{Bourbaki}.  
 For each pair of indices $(i,j)\in \mathcal{I}\times \mathcal{I}$,
 let $\zeta_{ij}$ and 
 $\phi_{ij}$ be defined so that the following diagram commutes:
$$\begin{CD}
C_0({\cal O}_X^{\oplus 3}|U_{ij}, q_i|U_{ij}) 
@>{\psi_{b_i}|U_{ij}}>{\cong}> \Lambda^{ev}({\cal O}_X^{\oplus 3}|U_{ij})
@>{\zeta_i|U_{ij}}>{\cong}> A|U_{ij}\\
@V{\phi_{ij}}V{\cong}V @V{\zeta_{ij}}V{\cong}V @V{=}VV\\
C_0({\cal O}_X^{\oplus 3}|U_{ij}, q_j|U_{ij}) 
@>{\cong}>{\psi_{b_j}|U_{ij}}> \Lambda^{ev}({\cal O}_X^{\oplus 3}|U_{ij})
@>{\cong}>{\zeta_j|U_{ij}}> A|U_{ij}
\end{CD}$$
The above diagram means that the algebras $A_i$ glue along 
 $U_{ij}:=U_i\cap U_j$ via $\zeta_{ij}$ to give (an algebra bundle 
 isomorphic to) $A$, and in the same vein, the even Clifford algebras 
 $C_0({\cal O}_X^{\oplus 3}|U_i, q_i)$ glue along the $U_{ij}$ via 
 $\phi_{ij}$ to give $A$ as well.  
Now consider the similarity $g_{l_{ij}}^{\phi_{ij}}=s_{-1}^{+}(\phi_{ij}):({\cal O}_X^{\oplus 3}|U_{ij}, q_i|U_{ij})\thinspace\cong_{l_{ij}}\thinspace({\cal O}_X^{\oplus 3}|U_{ij}, q_j|U_{ij})$ with multiplier $l_{ij}:=\hbox{\rm det}(\phi_{ij})^{-1}$ given by (c), 
 Theorem \ref{lifting-of-isomorphisms}. Since $s_{-1}^{+}$ 
 is  multiplicative, and since $\phi_{ij}$ satisfy 
 the cocycle condition, it follows that $s_{-1}^{+}(\phi_{ij})$ also 
 satisfy the cocycle condition and therefore glue the
 ${\cal O}_X^{\oplus 3}|U_i$
 along the $U_{ij}$ to give a rank 3 vector bundle $V$ on $X.$ While 
 the $q_i$ do not glue to give an ${\cal O}_X$-valued quadratic form on 
 $V$, the facts that the multipliers $\{l_{ij}\}$ form a cocycle for 
 $I:=\hbox{\rm det}^{-1}(A)$ and that 
 $s^+_{(-1)}$ is a section together imply, taking into account the 
 uniqueness in (1), Prop.\ref{simil-induces-iso-of-even-cliff},
  that actually the $q_i$ glue to give an 
 $I$-valued quadratic form $q$ on $V$ and 
 that $\covqi\cong A.$ It was verified in page 28, \cite{tevb-paper2} that 
  ${(\phi_{ij})}_{\Lambda^2}=\hbox{\rm det}(\phi_{ij}).
 \Lambda^2(g_{l_{ij}}^{\phi_{ij}}).$ This immediately implies part (1) of
 assertion (a) of Theorem \ref{structure-of-specialisation}, from which 
 parts (2)---(4) can be deduced using the standard properties of the
  determinant  
 and the perfect 
 pairings between suitable exterior powers of a bundle.

\paragraph{Proofs of assertions in (b),
 Theorem \ref{structure-of-specialisation}.}  We first prove (b1). 
Let $A$ be a given specialisation, and let $A\cong\covqi$ as in part 
 (a) of Theorem \ref{structure-of-specialisation} with
 $I=\hbox{\rm det}^{-1}(A).$ By the injectivity part of 
 Theorem \ref{bijectivity}, we have 
$\covqi\cong A\cong C_0(V',q',{\cal O}_X)$ iff 
there exists a twisted discriminant bundle $(L,h,J)$ and an isomorphism 
 $\vqi\cong (V',q',{\cal O}_X)\otimes (L,h,J).$ The latter implies that 
 $I\cong J\cong L^2$ and hence $\hbox{\rm det}(A)\in 2.\hbox{\rm Pic}(X).$ 
 On the other hand, if this last condition holds, we could take 
 for $L$ a square root of $J:=I^{-1}$, alongwith an isomorphism
 $h:L^2\cong J$ and we would have by Prop.\ref{isom-covq-covlqh}
 an algebra isomorphism 
$\gamma_{(L,h,J)}:C_0(V\otimes L, q\otimes h, {\cal O}_X)\equiv
C_0\left(\vqi\otimes(L,h,J)\right)\cong C_0(V,q,I)\cong A.$

For the proof of (b2), 
suppose that the line subbundle ${\cal O}_X.1_A\hookrightarrow A$ is a
 direct summand of $A.$  We
 may choose a splitting $A\cong {\cal O}_X.1_A\oplus 
 (A/{\cal O}_X.1_A).$ Using assertion (1) of (a),
 Theorem \ref{structure-of-specialisation}, we see that there exists a
 rank 3 vector bundle $V$ on $X$ such that
 $A\cong {\cal O}_X.1_A\oplus (A/{\cal O}_X.1_A)\cong {\cal O}_X.1_A\oplus 
 (\Lambda^2(V)\otimes I^{-1})\cong {\cal O}_X.1\oplus \Lambda^2(V)\otimes 
 I^{-1}=:W$ where $I:=\hbox{\rm det}^{-1}(A)$ and   
 the last isomorphism is chosen so as 
 to map ${\cal O}_X.1_A$ isomorphically onto ${\cal O}_X.1.$
 Therefore if $(W,w):=({\cal O}_X.1\oplus\Lambda^2(V)\otimes I^{-1},1)$,
 then  by the above identification $A$ induces an element of $\spazuWw(X)$, 
 and since $\Upsilon:\hbox{\rm Bil}_{(V,I)}
\cong \spazuWw$ is an $X$-isomorphism 
 by (1), Theorem \ref{bilinear-forms-as-specialisations}, it follows that 
 there exists an $I$-valued 
  global quadratic form $q=q_b$ induced from an $I$-valued global
  bilinear form 
  $b$ on $V$ such that the algebra structure $\Upsilon(b)\cong A.$ 
  (We recall that 
 $\Upsilon(b)$ is the algebra structure induced from the 
 linear isomorphism  $\psi_b:C_0(V,q=q_b, I)\cong {\cal O}_X.1\oplus
 \Lambda^2(V)\otimes I^{-1}=W$ of 
 (2d), Theorem \ref{Bourbaki}, which preserves 1 by (2a) of the same 
 Theorem). The proof of (b3) follows from a combination of 
 those of (b1) and (b2). {\bf Q.E.D., Theorem \ref{structure-of-specialisation} and 
 surjectivity part of Theorem \ref{bijectivity}.}

\begin{flushleft}
{\large \bf Acknowledgements}
\end{flushleft}
\addcontentsline{toc}{section}{Acknowledgements}
 The author gratefully acknowledges the Postdoctoral Fellowship 
 (September 2003--August 2004) of the {\em Graduiertenkolleg Gruppen
 und Geometrie} 
 under support from the {\em Deutschen Forschungsgemeinschaft}  and the 
 State of Niedersachsen at the Mathematisches Institut G\"ottingen 
 where this paper was written. The author also thanks the 
 National Board for Higher Mathematics, Department of Atomic Energy, 
 Government of India, for its Postdoctoral Fellowship (August 2002--August 2003) at the 
 Chennai Mathematical Institute, Chennai, India, during the latter half of 
 which certain special cases of the results of this work were obtained. 
 \vspace*{-4mm}

\end{document}